\documentclass[12pt]{amsart}

\usepackage{amsfonts}
\usepackage{latexsym}
\usepackage{amssymb}
\usepackage{amsthm}
\usepackage{amsmath,amssymb,amsfonts,latexsym,graphicx}
\usepackage{amsmath}
\usepackage{graphicx}
\usepackage{psfrag}
\usepackage{subfigure}
\usepackage{color}

\numberwithin{equation}{section} \numberwithin{figure}{section}
\def\n{\noindent}

\def\d{\displaystyle}

\def\bga{\begin{array}}
\def\eda{\end{array}}

\def\De{\Delta}
\def\pt{\partial}

\def\gm{\gamma}

\def\la{\lambda}

\def\sg{\sigma}

\def\al{\alpha}

\def\rw{\rightarrow}
\def\iy{\infty}

\def\d{\displaystyle}
\def\dfr#1#2{\displaystyle{\frac{#1}{#2}}}

 \theoremstyle{definition}
 
 \theoremstyle{remark}

\def\bp{\mathbf{p}}

\def\bF{\mathbf{F}}

\def\bI{\mathbf{I}}

\def\bP{\mathbf{P}}

\def\bU{\mathbf{U}}

\def\beqn{\begin{equation}}
\def\eeqn{\end{equation}}

\title{Stiffened gas approximation and \\ GRP resolution for fluid flows of real materials }

\thanks{   }

\author{YUE WANG}
\address{Yue Wang:
 Laboratory of Computational Physics, Institute of Applied Physics and Computational Mathematics, Beijing, 100088, China  }
 \email{Yue Wang: wang\_yue@iapcm.ac.cn}

\author{JIEQUAN LI}
\address{Jiequan Li: Laboratory of Computational Physics, Institute of Applied Physics and Computational Mathematics, Beijing, 100088, China; and
HEDPS, Center for Applied Physics and Technology, and College of Engineering, Peking University, Beijing 100871, China}
\email{
Jiequan Li:  li\_jiequan@iapcm.ac.cn
}

\begin{document}



\begin{abstract}
The equation of state (EOS) embodies thermodynamic  properties of compressible fluid materials and  usually has very complicated forms in real engineering applications, subject to the physical requirements of thermodynamics. The complexity of EOS in form gives rise to the difficulty in analyzing relevant wave patterns.  Concerning the design of  numerical algorithms, the complex EOS causes the inefficiency of Riemann solvers and even  the loss of robustness, which hampers the  development of   Godunov-type  numerical schemes.
 In this paper, a strategy of local stiffened gas approximation is proposed for real materials.  The stiffened gas EOS is used to approximate general EOS locally at each interface of computational control volumes so that the Riemann solver can be significantly simplified.  In the meantime, the generalized Riemann problem (GRP) solver is adopted not only for high resolution purpose but effective reflection of the local thermodynamics as well. The resulting scheme is demonstrated to be  efficient and robust and numerical examples display the excellent performance of such an approximation.

\end{abstract}
\maketitle


\pagestyle{myheadings} \thispagestyle{plain} \markboth{
  YUE WANG and JIEQUAN LI} { Stiffened gas approximation and GRP resolution for fluid flow of real materials }

\section{Introduction}\label{sec1}

The dynamical evolution of compressible fluid flows is determined by the conservation laws of mass, momentum, and energy.
To obtain a complete mathematical description of fluids, the conservation laws must be supplemented with the equation of state (EOS) besides constitutive relations characterizing the material properties. In this context,  the fundamental Gibbs relation
\begin{align}
TdS=de+pd\tau \label{Gibbs}
\end{align}
connects all thermodynamic variables, the temperature $T$, the entropy $S$, the internal energy $e$, the pressure $p$, the density $\rho$  and the specific volume $\tau=1/\rho$.  This relation also  shows that   just two of them are independent.  We choose $\rho$ and $S$ as the independent variables and express the EOS in a complete form, e.g.,
\begin{align}
p=p(\rho,S).
\label{rho-s}
\end{align}
In practice, the EOS is often expressed in terms of another pair of independent variables
\begin{align}
p=p(\rho,e)
\label{rho-e}
\end{align}
such as the Mie-Gr{\"{u}}neisen EOS \cite{Plohr}.  Our familiar examples are the EOS for polytropic gases
\begin{equation}
p=(\gamma-1)\rho e, \ \ \ \gamma >1,
\end{equation}
and stiffened gases
\begin{equation}
p=(\gamma-1)\rho e - \gamma p_\infty,\label{stiff-EOS}
\end{equation}
where $p_\infty$ is a constant.  Note that the $(\rho, e)$-formula of EOS can be expressed  in terms of $(\rho, S)$ by using the Gibbs relation \eqref{Gibbs}  \cite{Plohr}.  In this paper we consider the following family of EOSs.

\n {\bf Assumption.}  {\em The pressure $p$ has a linear relation with the internal energy $e$,
\begin{align}
p=p(\rho, e) =\kappa(\rho)e +\chi(\rho), \label{EOS}
\end{align}
where $\kappa(\rho)$ and $\chi(\rho)$ are two given functions of  density $\rho$.
The EOS describing the fluid is convex \cite{Plohr}, i.e.,  the fundamental derivative
\begin{align}
\mathcal{G}:=-\frac 12\tau[(\pt^2p/\pt \tau^2)_S/(\pt p/\pt\tau)_S] \label{fundamental_derivative}
\end{align}
 is strictly positive.
 }

The family of EOS \eqref{EOS} has extensive applications as shown in numerical examples in Section \ref{sec:num}.
With this EOS, the Euler system of  conservation laws consisting  of mass, momentum and energy  is closed and can be used to describe the dynamics of compressible fluid flows, for which the Riemann problem plays a fundamental role not only in understanding  the wave patterns but also in  designing numerical algorithms.
There are abundant literatures about how to solve the Riemann problem \cite{banks,Godunov,Plohr,Quar-2003,Thompson,Toro}.
Many of the associated relations can be evaluated in closed form for the polytropic gases. However, several of these relations must be solved numerically for general EOSs.
Kamm \cite{Kamm} summarizes an algorithm for solving the Riemann problem of compressible flow with general convex EOSs.
The approach for the  Jones-Wilkins-Lee (JWL) EOS is coded in the ExactPack package to derive exact solutions for code verification. Since there involve an appreciable number of EOS calls  in the procedure of  finding the zeroes of the nonlinear equations and numerically solving ODEs for rarefaction wavces, this approach is prohibitively inefficient as the ``exact" Riemann solver for a general compressible flow algorithm.
As far as approximated Riemann solvers are concerned, the oversimplification may destroy the thermodynamic property that the original EOS embodies.
Hence, the complexity of EOS hampers the practical application of associated Riemann solvers and the development of corresponding Godunov-type numerical schemes.

The EOS embodies  the subtle thermodynamic properties of compressible fluids. The interaction of   kinematics  and thermodynamics determines  the whole dynamical process. The more significant  the thermodynamic effect is  in the fluid dynamics, the more important role the EOS plays. For  numerical purpose in practice,   proper approximation has  to be made  so that numerical fluxes can be constructed effectively. However, we have to bear in mind that any approximation should inherit and preserve the inherent thermodynamic property of the EOS.   To this end  we propose a new strategy, {\bf a Stiffened Gas Approximation} of EOS for real materials,  in order to suit for the design of high resolution Godunov-type schemes.
 Considering EOS as a part of the dynamical system of fluids, the stiffened gas approximation means that at each local background state $\rho_0$, the EOS \eqref{EOS} is approximated as
\begin{align}
&p=\tilde{p}(\rho,e)=(\gm(\rho_0)-1)\rho e-\gm(\rho_0)p_{\infty}(\rho_0),\label{approx-stiffen}\\
&\gamma(\rho_0)=1+\kappa(\rho_0)/\rho_0,\ p_{\infty}(\rho_0)=-\chi(\rho_0)/\gamma(\rho_0).\label{approx-stiffen1}
\end{align}
As such an approximation is applied to the Riemann problem for compressible  fluid flows of real materials, a new {\em two-material Riemann problem} is formulated. As far as this strategy is applied for Godunov-type numerical schemes, the stiffened gas approximation is made over each control volume and the two-material Riemann problem is solved at each cell boundary with approximate stiffened gas EOSs. It turns out that the local Riemann solver at each cell boundary has the same simplicity as that for polytropic gases.
It is noted that  a similar idea can be found in \cite{Saurel} to derive the discretization of the non-conservative equation in a two-phase flow model. However, significance of the stiffened gas approximation and its application on numerical fluxes have not been discussed yet.

In order to manifest the thermodynamic properties  numerically that EOS embodies,  we adopt the generalized Riemann problem (GRP) solver locally \cite{Ben-Artzi-1984,Ben-Artzi-2003,Ben2006}.  As we discussed in \cite{Li-Wang},   the GRP solver  not only plays  the role of accuracy in the design of high resolution Godunov-type schemes, but also provides an approach including the thermodynamic effect in the resulting scheme.  The passage is taken through precisely measuring the entropy variation (by the Gibbs relation)  and its interaction with kinematic and other thermodynamic variables.   Such consideration provides an evident and immediate performance  when thermodynamic effect is severe such as in the LeBlanc problem (the  problem for large ratios of density or pressure \cite{Tang-Liu, Li-Wang}).
 As an immediate product of the stiffened gas approximation to the Riemann problem and the generalized Riemann problem (GRP),   a new type of approximate Godunov type scheme is designed for general EOSs \eqref{EOS}. The main advantages using the stiffened gas approximation are listed as follows.

\begin{enumerate}
  \item[(i)]  Difficulty of the exact Riemann solver with general EOS is avoided, thus the resulting scheme can be more efficient and robust.
  \item[(ii)]  The thermodynamic relations are derived explicitly  and  manifested in the design of the approximate GRP solver to guarantee resulting numerical solutions obey the same or at least an approximate analogue.
  \item[(iii)] The stiffened gas approximation is made just locally at each interface of control volumes. Therefore  the property of EOS is maintained as much as possible (even in form).
\end{enumerate}

This paper is organized by six sections.
Besides the introduction here, we will show the stiffened gas approximation of the EOS and formulate an approximate two-material Riemann problem with stiffened gases instead of solving a single-material Riemann problem with general EOS in Section 2.
The stiffened gas approximation of the generalized Riemann problem is presented in Section \ref{sec-GRP}  and  the issue of thermodynamics is discussed during the approximation.
Numerical schemes are summarized in Section 4 and  numerical examples are presented in Section 5 to demonstrate the performance of the resulting scheme. Concluding remarks are given in Section 6.

\newpage

\section{stiffened gas approximation of the Riemann problem}\label{sec-RP}

The Euler equations for compressible fluids are written in the form
\begin{equation}
\begin{array}{c}
\dfr{\pt \bU}{\pt t}  +\dfr{\pt \bF(\bU)}{\pt x}=0,\\[3mm]
\bU =(\rho, \rho u,   \rho E)^\top,  \ \ \ \bF(\bU) = (\rho u, \rho u^2 +p,  u(\rho E +p))^\top,
\end{array}
\label{Euler-1D}
\end{equation}
where  the total energy $E=u^2/2 +e$   consists of the internal energy $e$ and the kinematic energy $u^2/2$.
 The system is closed by an equation of state using two independent thermodynamical variables.
Here we consider the convex equation of state in the form of \eqref{EOS}
\begin{align}
p=p(\rho,e)=\kappa(\rho)e +\chi(\rho).\notag
\end{align}
For instance, the stiffened gas takes $\kappa(\rho) = (\gamma-1)\rho$ and $\chi(\rho)=-\gamma p_{\infty}$ where $\gamma$ and $p_{\infty}$ are constant parameters \cite{Godunov1979}.
It can be reduced to the polytropic gas with $p_{\infty}=0$.

The Riemann problem for \eqref{Euler-1D} is subject to the initial data
\begin{equation}
\begin{array}{l}
 \bU(x,0) =\left\{
 \begin{array}{ll}
 \bU_L, & x<0,\\
 \bU_R, & x>0.
 \end{array}
 \right.
 \end{array}
 \label{RP}
\end{equation}
As we solve the Riemann problem \eqref{RP} numerically with the original complex EOS \eqref{EOS},   computational complexity is unavoidable, as  mentioned in \cite{Kamm} for the exact Riemann solver:
\begin{enumerate}
  \item[(i)]  numerical integration and interpolation for the rarefaction waves;
    \item[(ii)]  nonlinear equations with the intricate relations of shocks and rarefaction waves;
  \item[(iii)]reasonable ``starting" estimates of the star-state pressure;
  \item[(iv)]  many EOS calls during the procedure for solving nonlinear equations.
\end{enumerate}
Due to the above reasons, the exact Riemann solver for complex EOSs is prohibitively inefficient  for a general
compressible flow algorithm.

In order to make the Riemann solver more efficient and robust as well as maintain details of wave interactions of the Riemann solution, we make the following stiffened gas approximation and formulate
{\bf an approximate two-material stiffened gas Riemann problem} for \eqref{Euler-1D} and \eqref{RP}.
Specifically, we have two cases.

\begin{enumerate}
\item[(i)] {\bf Case 1: $u_L=u_R$ and $p_L=p_R$}. There is no need making any approximation. The Riemann solution is
\begin{equation}
(\rho,u,p)(x,t) =\left\{
\begin{array}{ll}
(\rho_L,u_L,p_L), \ \ \ &x<u_Lt=u_R t,\\
(\rho_R, u_R, p_R), & x>u_Lt=u_R t.
\end{array}
\right.
\end{equation}
It turns out that the internal energy is defined through the EOS
\begin{equation}
\begin{array}{ll}
p_L =\kappa(\rho_L) e(x,t) +\chi(\rho_L),  \ \ \ &x<u_Lt=u_R t,\\
p_R =\kappa(\rho_R) e(x,t) +\chi(\rho_R), & x>u_Lt=u_R t.
\end{array}
\end{equation}

\item[(ii)]  {\bf Case 2: $u_L\neq u_R$ or  $p_L\neq p_R$}.  The stiffened gas approximation is applied for the EOS locally. Then a two-material Riemann problem is formulated. The EOS of the fluid in the left-hand side of the material interface is approximated as
\begin{equation}
   p(\rho,e)=
 (\gm(\rho_L)-1)\rho e-\gm(\rho_L)p_{\infty}(\rho_L);
 \label{stiff-L}
\end{equation}
while the EOS in the right-hand side is
\begin{equation}
   p(\rho,e)=
 (\gm(\rho_R)-1)\rho e-\gm(\rho_R)p_{\infty}(\rho_R).
 \label{stiff-R}
\end{equation}
Here $\gamma(\rho_i)=1+\kappa(\rho_i)/\rho_i,\ p_{\infty}(\rho_i)=-\chi(\rho_i)/\gamma(\rho_i),\ i=L,R$.
Note that the material interface is determined together with the solution and not known a priori, the same as the contact discontinuity for a single-material Riemann problem.

\end{enumerate}

The outcome of the stiffened gas approximation is evident.
Since the stiffened gas is a typical type of the general EOS \eqref{EOS}, the Riemann solution of \eqref{Euler-1D} and \eqref{RP} is  self-similar, consisting of a contact discontinuity and a combination of (expansive) rarefaction waves and (compressive) shocks. This approximate two-material stiffened gas Riemann problem does not give rise to  any new difficulty. The corresponding Riemann solver is analogous to that
 for  polytropic gases  \cite{Ben-Artzi-2003, Godunov,Godunov1979,Gottlieb,Toro}.   Let us take  a Rarefaction-Contact-Shock structure as an example. Then Riemann invariants and a isentropic relation provide a key for the expansion of  the  rarefaction wave in the left
\begin{equation}
\begin{array}{l}
u_L+\frac{2c_L}{\gm(\rho_L)-1}=u_*+\frac{2c_{*L}}{\gm(\rho_L)-1},\\
(p_L+p_{\infty}(\rho_L))\rho_L^{-\gm(\rho_L)} = (p_*+p_{\infty}(\rho_L))\rho_{*L}^{-\gm(\rho_L)},
\end{array}
\label{RW}
\end{equation}
where $c$ is sound speed, given by $c^2 =\gm (p+p_\iy)/\rho$ .
The Rankine-Hugoniot (RH) relations provide the basis for resolution of the shock wave in the right,
\begin{equation}
\begin{array}{l}
 \ [\rho u] = \sigma[\rho],\\
\ [\rho u^2+p] = \sigma[\rho u],\\
 \ [u(\rho E+p)] = \sigma[\rho E],
\end{array}
\label{SW}
\end{equation}
where $[\cdot] = (\cdot)_R-(\cdot)_{*R}$ and $\sigma$ is the shock speed. We combine \eqref{RW} and \eqref{SW} to obtain
\begin{align}
&u_* = u_L-\frac{2c_L}{\gm(\rho_L)-1}\left[(\frac{p_*+p_{\infty}(\rho_L)}{p_L+p_{\infty}(\rho_L)})^{\frac{\gm(\rho_L)-1}{2\gm(\rho_L)}}-1\right]=: u_L-f_L(p_*,\rho_L,p_L),\notag\\
&u_* = u_R+(p_*-p_R)\left[\frac{\frac{2}{\gm(\rho_R)+1}\frac{1}{\rho_R}}{p_*+\frac{\gm(\rho_R)-1}{\gm(\rho_R)+1}p_R+\frac{2\gm(\rho_R)}{\gm(\rho_R)+1}p_{\infty}(\rho_R)}\right]^{\frac 12}=: u_R+f_R(p_*,\rho_R,p_R).\notag
\end{align}
The equality of the pressure and velocity across the material interface provides the bridge to solve the Riemann problem exactly
\begin{align}
\Delta u + f_R(p_*,\rho_R,p_R)+f_L(p_*,\rho_L,p_L)=0,\ \Delta u = u_R-u_L.\notag
\end{align}
So any nonlinear iteration scheme for this algebraic equation is the same as that for polytropic gases in the Riemann solver \cite{Toro}.

We  close this section by concluding the gains from the stiffened gas approximation for the Riemann problem.
\begin{enumerate}
  \item For the Riemann problem with general EOS,  a secant method is needed since the formulations can not be written explicitly. The Newton-Raphson iterative procedure can be applied for the stiffened gas approximation so that the iteration will be more efficient.
  \item For the exact Riemann solver of complex EOS, except the iteration for the star state $p_*$ and $u_*$, nonlinear equations are also solved inside the rarefaction waves and shocks. However, only one nonlinear equation is needed for the approximation case. As a result, the EOS calls are substantially reduced.
  \item Since there is  the explicit relation for  rarefaction waves, numerical integration and interpolation are avoided.
  \item The Riemann problem is still ``precisely" solved since all of the information about the wave configuration is known.
\end{enumerate}

\section{Stiffened gas approximation of the generalized Riemann problem and thermodynamic effects}\label{sec-GRP}

In this section, we will discuss the generalized Riemann problem (GRP) for \eqref{Euler-1D} subject to the discontinuous initial data
\begin{equation}
\begin{array}{l}
 \bU(x,0) =\left\{
 \begin{array}{ll}
 \bP_L(x), & x<0,\\
 \bP_R(x), & x>0,
 \end{array}
 \right.
 \end{array}
 \label{GRP}
\end{equation}
where $\bP_L(x)$ and $\bP_R(x)$ are  smooth functions, but $\bP_L(0-0)\neq \bP_R(0+0)$ generally speaking.   The GRP solver serves to compute instantaneous values
\begin{equation}
\bU_* = \lim_{t\rw 0^+} \bU(0,t), \ \ \ \left(\dfr{\pt \bU}{\pt t}\right)_* = \lim_{t\rw 0^+} \dfr{\pt \bU}{\pt t}(0,t).
\label{GRP-value}
\end{equation}
Denote $\bU_L = \bP_L(0-0)$, $\bU_R=\bP(0+0)$, and $\bU_L'=\bP'(0-0)$,$\bU'_R=\bP'(0+0)$.
 Our previous studies \cite{Ben-Artzi-1984,Ben-Artzi-2003,Ben2006} show that the instantaneous values \eqref{GRP-value} only depend on the limiting values $(\bU_L, \bU_R)$ and $(\bP'_L(0-0), \bP_R'(0+0))$.  Therefore we might as well  discuss the ``linear GRP" for \eqref{Euler-1D} directly  with the piecewise linear distribution
\begin{equation}
\bU(x,0)=\left\{
\begin{array}{ll}
\bU_L +\bU_L'x, \ \ \ &x<0,\\
\bU_R+\bU_R'x, \ \ \ &x>0,
\end{array}
\right.
\label{data-GRP}
\end{equation}
and compute the corresponding values \eqref{GRP-value}.
Such a discussion is not limited to the theoretical significance. It is  the basis of high resolution numerical schemes \cite{Ben-Artzi-1984,Ben-Artzi-2003,Ben2006,Ben2007}. More importantly, the thermodynamic effect is included in numerical fluxes through the GRP solver \cite{Li-Wang}, which is a main issue that we focus on now.

\begin{figure}[!htb]
\centering \subfigure[Wave pattern for the GRP  ]{
\psfrag{0}{$0$}
\psfrag{x}{$x$}\psfrag{t}{$t$}\psfrag{betal}{\small
$\beta=\beta_L$}\psfrag{ul}{$\bU_L$}\psfrag{al2}{$\bar{\bar\al}$}
\psfrag{shock}{shock}\psfrag{rarefaction}{\small rarefaction}\psfrag{uminus}{$\bU_L(x,t)$}\psfrag{uplus}{$\bU_R(x,t)$}
\psfrag{u1}{$\bU_{1*}$}\psfrag{u2}{$\bU_{2*}$}\psfrag{al=albar}{$\al=\bar\al$}\psfrag{xbar}{$\bar
\al$ }\psfrag{al=2albar}{$\al=\bar{\bar \al}$}\psfrag{ur}{$\bU_R$}
\psfrag{contact}{contact}\psfrag{betastar}{$\beta=\beta_*$}
\psfrag{ul}{$\bU_L$}\psfrag{ur}{$\bU_R$}\psfrag{ustar}{$\bU_*$}
\includegraphics[clip, width=7.6cm]{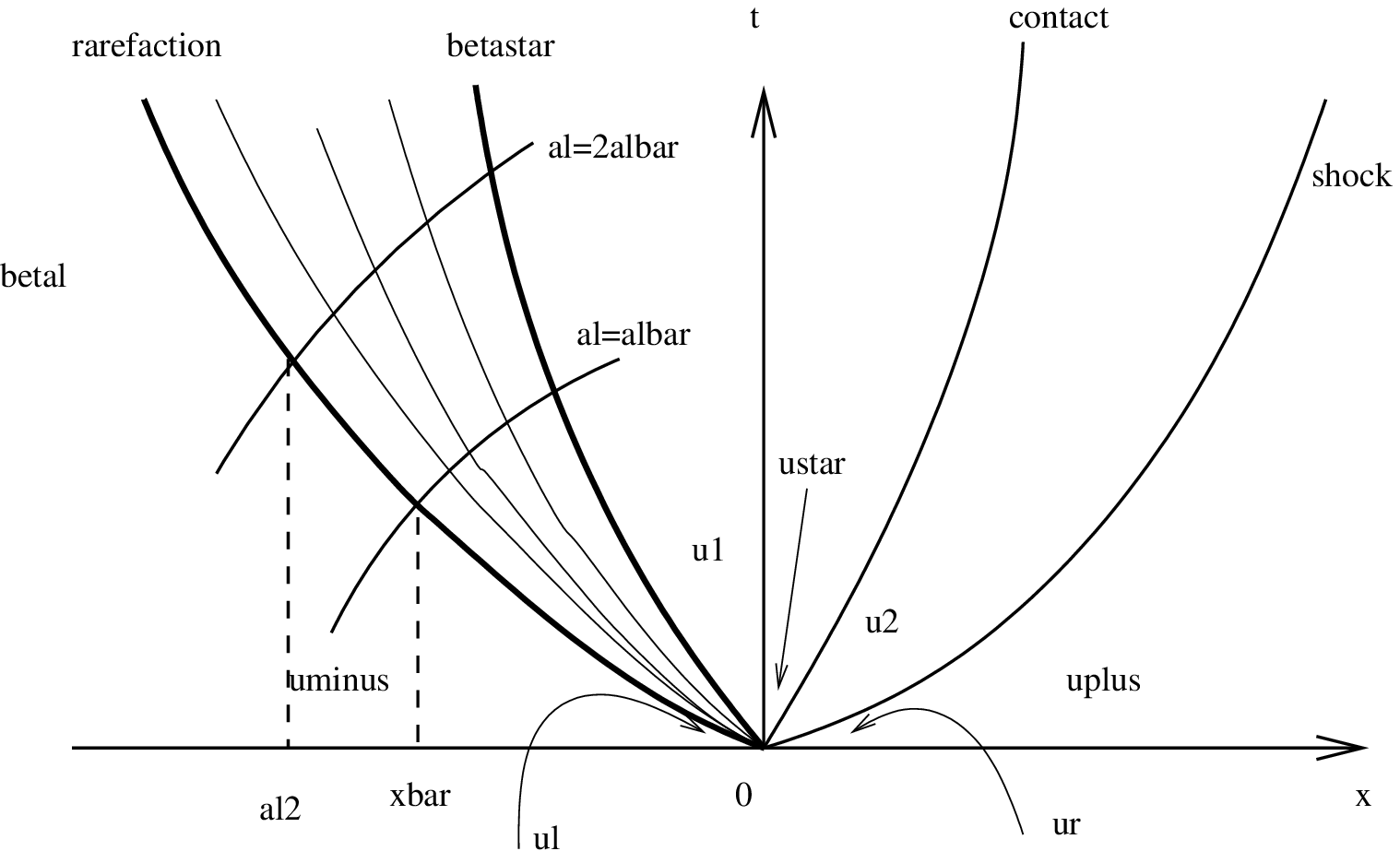} }
\hspace{0.5cm}
 \subfigure[Wave pattern for the associated RP]{
\psfrag{x}{$x$}\psfrag{t}{$t$}\psfrag{betal}{\small
$\beta=\beta_L$}
\psfrag{shock}{shock}\psfrag{rarefaction}{rarefaction}\psfrag{uminus}{$\bU_L$}\psfrag{uplus}{$\bU_R$}
\psfrag{u1}{$\bU_{1*}$}\psfrag{u2}{$\bU_{2*}$}\psfrag{al=albar}{$\al=\bar\al$}\psfrag{xbar}{$\bar
\al$ }\psfrag{al=2albar}{$\al=\bar{\bar \al}$}
\psfrag{contact}{contact}\psfrag{betastar}{$\beta=\beta_*$}\psfrag{ustar}{$\bU_*$}
\psfrag{0}{$0$}\psfrag{al2}{$\bar{\bar \al}$}
\includegraphics[width=8.5cm]{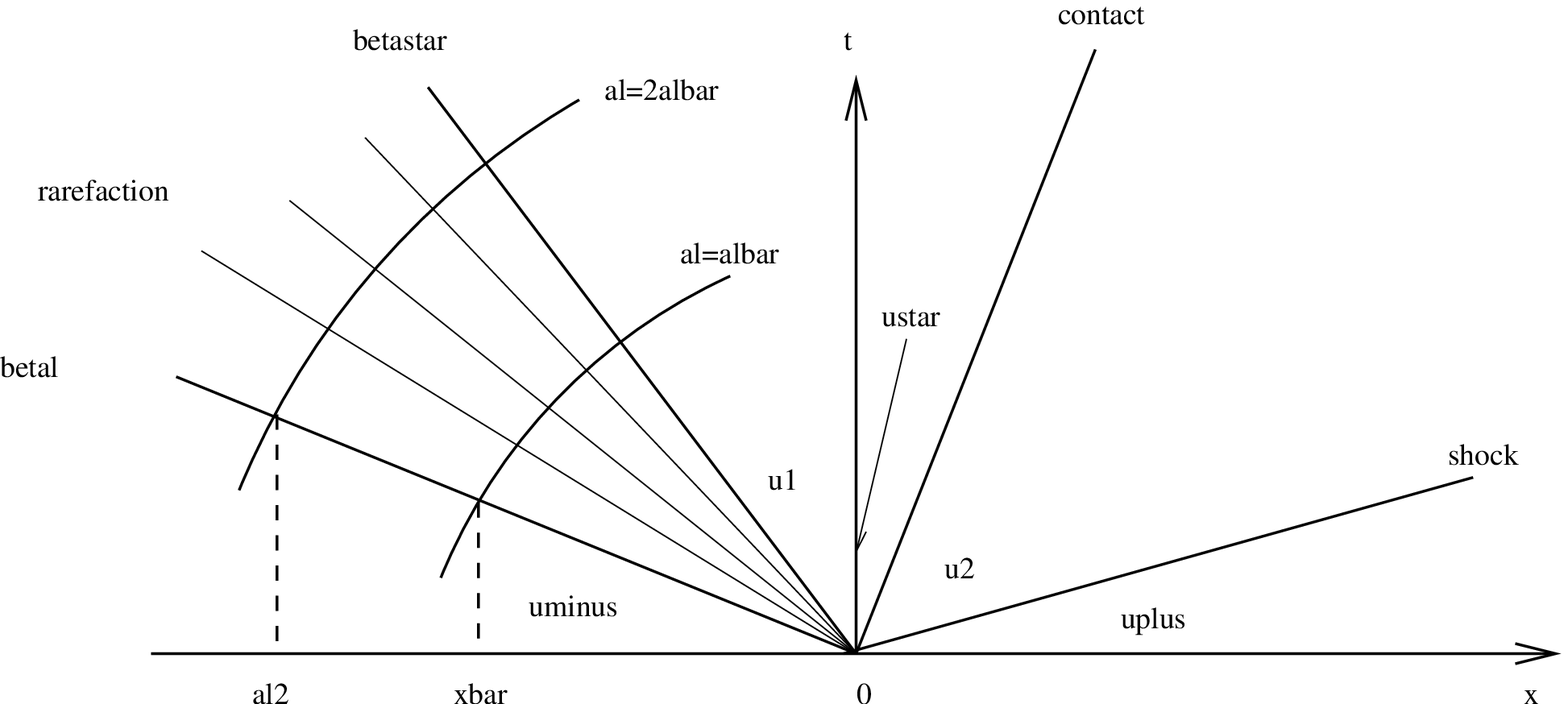}}
\caption[small]{Typical wave pattern for the generalized Riemann problem and its connection with the associated Riemann problem.}
\label{Fig-wave}
\end{figure}

As justified in \cite{Ben-Artzi-2003, Ben2007}, the resolution of GRP is based on the connection with the associated Riemann problem, i.e., the Riemann problem for \eqref{Euler-1D}   subject to the initial data $(\bU_L, \bU_R)$.
This associated Riemann solution is denoted as $\bU^A(x/t) =R^A(x/t;\bU_L,\bU_R)$, and the solution of GRP \eqref{Euler-1D}-\eqref{GRP} is $\bU(x,t)$. Then we have
\begin{equation}
\lim_{t\rw 0^+}\bU(\la t,t )=R^A(\la;\bU_L,\bU_R), \ \  \ \ \la>0.
\end{equation}
 In Figure \ref{Fig-wave}, we display a wave configuration of rarefaction-contact-shock for GRP and its connection with the associated Riemann solution.
With such a connection, we see that the local wave configuration is determined by the associated Riemann solution. Hence the local two-material Riemann problem,  i.e., the stiffened gas approximation with     the background states $\bU_L$ and $\bU_R$, is formulated to approximate  the local wave configuration. The value $\bU_*=R^A(0;\bU_L, \bU_R)$ is computed using the two-material Riemann solver that we proposed in the last section.

%

Based on this,  the GRP can be  used to compute the instantaneous derivative value in \eqref{GRP-value} and resolved  for the approximate stiffened gas generalized Riemann problem.  We adopt the direct Eulerian version in \cite{Ben2006}.
With reference to Figure \ref{Fig-wave},  the left-going rarefaction wave and the right-going shock should be resolved, respectively.  Details are put in Appendix.  As we discussed in \cite{Li-Wang}, the GRP solver is not only used to increase accuracy of order for  high resolution Godunov-type (e.g. GRP)  schemes,  but also reflect the thermodynamic effect. Let's comment on this issue from
the entropy variation in terms of the EOS.

 In the case of EOS \eqref{EOS}, we have
\begin{equation}
dp= \dfr{\pt p(\rho, e)}{\pt \rho} d\rho +\dfr{\pt p(\rho, e)}{\pt e} d e.
\end{equation}
The Gibbs relation implies
\begin{equation}
\dfr{\pt p(\rho, e)}{\pt e} TdS = dp -\left( \dfr{\pt p(\rho, e)}{\pt e}  p/\rho^2+\dfr{\pt p(\rho, e)}{\pt \rho} \right) d\rho.
\end{equation}
Specified to the EOS  \eqref{EOS}, we have
\begin{equation}
dp =( \kappa'(\rho) e+\chi'(\rho)  )d\rho+ \kappa (\rho)de.
\end{equation}
Hereafter the prime $'$ refers to the corresponding derivative with respect to its independent variable, and the explanation is omitted when no confusion is caused.
We continue to use the Gibbs relation \eqref{Gibbs} and obtain all differential relations, e.g.,
\begin{equation}
\kappa(\rho) TdS =dp - ( \kappa'(\rho) e+\chi'(\rho)  +\kappa(\rho) p/\rho^2 )d\rho.
\end{equation}
Hence the entropy variation can be expressed in terms of the variations of primitive variables $(\rho, p)$ . In particular,  the initial entropy variations are measured as
\begin{equation}
\kappa(\rho_i) T_iS_i' =p_i' - ( \kappa'(\rho_i) e_i+\chi'(\rho_i)  +\kappa(\rho_i) p_i/\rho_i^2 )\rho_i', \ \ \ i=L,R,
\label{ent-ini}
\end{equation}
and $e_i =e(p_i,\rho_i)$ by using \eqref{EOS}. Hence at the GRP level, there is no need making any  approximation about the entropy variation.
\vspace{0.2cm}

As the left-going rarefaction wave is considered, the initial entropy variation gives rise to the future variation. As in \cite{Li-Wang} (the notations are adopted there and also in Appendix), we find that there holds
\begin{equation}
\dfr{T (\pt S/\pt x)_*}{T_LS_{L}'}  =\theta_*^{\frac{1}{\mu_L^2}+1}, \ \ \theta_* =\frac{c_{*L}}{c_L}, \ \ \mu_L^2=\mu^2(\rho_L)=\frac{\gm(\rho_L)-1}{\gm(\rho_L)+1},  \label{s-al-poly}
\end{equation}
where $c$ is the local sound speed for approximate stiffened gases,
\begin{equation}
\ c^2 =\dfr{\gm (p+p_\iy)}{\rho}.
\end{equation}
The ratio \eqref{s-al-poly} shows the entropy evolution  across the (curved in GRP) non-isentropic rarefaction wave and reflects the thermodynamic effect, which further affects the variation of  other (kinematic  and thermodynamic) variables.

To see this point, we write the Euler equations in terms of so-called ``Riemann invariants" $(\phi, \psi,S)$,
\begin{equation}
\begin{array}{l}
\d \dfr{\pt\phi}{\pt t}+(u-c)\frac{\pt\phi}{\pt x}=T\dfr{\pt S}{\pt x},\\[3mm]
\d \frac{\pt\psi}{\pt t}+(u+c)\frac{\pt\psi}{\pt x}=T\dfr{\pt S}{\pt x},\\[3mm]
\dfr{\pt S}{\pt t}+u\dfr{\pt S}{\pt x} = 0,
\end{array}
\label{diagonal}
\end{equation}
where the Riemann invariants $(\phi,\psi)$  are
\begin{align}
\phi = u - \int^{\rho}\frac{c(\omega,S)}{\omega}d\omega, \ \ \psi = u + \int^{\rho}\frac{c(\omega,S)}{\omega}d\omega,\ \notag
\end{align}
with differential relations
\begin{equation*}
d\phi=du-\dfr{1}{\rho c} dp -\dfr 1c T dS, \ \ \ d\psi
=du+\dfr{1}{\rho c} dp +\dfr 1c TdS.
\end{equation*}
From \eqref{diagonal}, it is observed that unless the flow is isentropic or the initial distribution of data is uniform, the entropy variation ($T\pt S/\pt x\neq 0$) plays an important role  in the dynamics.  We still exemplify the wave configuration in Figure \ref{Fig-wave} and obtain
\begin{equation}
\left(\dfr{D u}{Dt}\right)_* + \dfr{1}{\rho_{*L} c_{*L}} \left(\dfr{D p}{Dt}\right)_* =d_L,  \ \ \dfr{D}{Dt} =\dfr{\pt }{\pt t} +u_*\dfr{\pt }{\pt x},
\label{rar-L}
\end{equation}
where  $d_L$ has an explicit expression
\begin{equation*}
d_L=\left[\dfr{1+\mu_L^2}{1+2\mu_L^2}
\theta_*^{1/(2\mu_L^2)}+\dfr{\mu_L^2}{1+2\mu_L^2}
\theta_*^{(1+\mu_L^2)/\mu_L^2}\right]T_LS'_L-c_L\theta_*^{1/(2\mu_L^2)}
\psi'_L.
\end{equation*}
From this expression, we see immediately how the initial entropy variation (the variation of initial distribution) takes effect on the variation of velocity $u$ and the pressure $p$ quantitatively, i.e.,  the influence of thermodynamics on kinematics.

 For the right-moving shock wave associated with $u+c$, the Rankine-Hugoniot relations \eqref{SW} hold and consist of two parts:  The Hugoniot condition
  \begin{equation}
  e_{*R}-e_R+(\tau_{*R}-\tau_R)\dfr{p_*+p_R}2=0,
 \end{equation}
 describes the relation of thermodynamic variables across the shock; and the kinematic-thermodynamic  jump relation
 \begin{equation}
 [\rho u]^2 =[\rho][\rho u^2+p],
 \end{equation}
 describes how the thermodynamic change affects the kinematic jump, and vice versa.  Here notations refer to \eqref{SW} and  Appendix. We track the shock trajectory to obtain
 \begin{equation}
a_R\left(\dfr{Du}{Dt}\right)_* +b_R\left(\dfr{D p}{Dt}\right)_* =d_R, \label{shock-R}
\end{equation}
where $a_R$, $b_R$ and $d_R$ are put in Appendix. We remark that there is no need approximating the EOS for the derivation of \eqref{shock-R}. Indeed,  the EOS \eqref{EOS} provides,
\begin{equation}
dp =\kappa(\rho) de + (\kappa'(\rho) e +\chi'(\rho))d\rho.
\end{equation}
The energy equation in \eqref{Euler-1D} can be written as
\begin{equation}
\dfr{\pt e}{\pt t} +u\dfr{\pt e}{\pt x} +\dfr{p}{\rho}\dfr{\pt u}{\pt x}=0.
\end{equation}
Therefore, as the Lax-Wendroff approach (the exchange of temporal and spatial derivatives)  is used for \eqref{shock-R}, the stiffened gas approximation is not necessary for the derivation of \eqref{shock-R}.  Of course, the adoption of stiffened gas approximation can simplify the calculation.

Thanks to the fact that the continuity of  $u$ and $p$ across the material interface implies the continuity of  the corresponding material derivatives, we combine \eqref{rar-L} and \eqref{shock-R} to  obtain $\left(\frac{D u}{Dt}\right)_*$, $\left(\frac{D p}{Dt}\right)_*$, and subsequently $\left(\frac{\pt u}{\pt t}\right)_*$, $\left(\frac{\pt p}{\pt t}\right)_*$.  The value $\left(\frac{\pt \rho}{\pt t}\right)_*$ follows from the EOS. Thus the instantaneous derivative value in \eqref{GRP-value} is  available.

As a summary, we comment on the thermodynamic effects by two points.
\begin{enumerate}
  \item The  Gibbs relation brings the initial thermodynamics into the entropy variation across the rarefaction waves. The ``kinematic-thermodynamic" variables $\psi$ and $\phi$ provide a bridge between kinematical and thermodynamic variables.
  \vspace{0.2cm}

  \item Using the Rankine-Hugoniot relations to track the singularity, the relations between kinematics and thermodynamics are built automatically across for shock waves.
\end{enumerate}

\section{Godunov type schemes with stiffened gas approximation}

Once the stiffened gas approximation for the general EOSs \eqref{EOS} is available,  the Godunov type scheme for \eqref{Euler-1D} is immediate. We just state the new scheme in one space dimension, and there is no difference at the point of the approximation of EOS for  the extension in multiple space dimensions.

Assume that the spatial computational domain $\Omega \in \mathbb{R}$ is covered completely by pairwise disjoint spatial elements
$I_j\in[x_{j-\frac 12},x_{j+\frac 12}]$, with $\De x_j=x_{j+\frac 12}-x_{j-\frac 12}$ and the cell average of $\bU$ within $I_j$ is defined at time $t=t_n$ as
\begin{align}
\bU_j^n=\frac{1}{\De x_j}\int_{I_j}\bU(x,t_n)dt.\notag
\end{align}
A typical Godunov type scheme for \eqref{Euler-1D} assumes a piecewise polynomial approximation at each time step $t=t_n$ over control volume $I_j$,
\begin{equation}
\bU(x,t_n) = \bp_j^{n}(x)\in \bP^k(I_j), \ \ \  x\in I_j,
\label{data}
\end{equation}
where $\bP^k(I_j)$ is the space of polynomials of degree $k$ on the cell $I_j$.
Let $\De t$ be the time step size and satisfy the usual CFL constraint.
The numerical solution is updated in two steps.
\begin{enumerate}
\item[(i)] {\em Average advancing.}  We advance the solution average  of \eqref{Euler-1D} and \eqref{data} to the next time step  $t_{n+1} =t_n +\De t$ according to the formula
\begin{align}
&\bU_j^{n+1} =\bU_j^n -\dfr{\De t}{\De x} [\bF_{j+\frac 12}^* - \bF_{j-\frac 12}^*],
\label{balance-scheme}
\end{align}
where $\bF_{j+\frac 12}^*$ is the numerical flux with desired accuracy at the cell boundary $x=x_{j+\frac 12 }$.

\vspace{0.2cm}

\item[(ii)] {\em Projection of data.}  We project the solution $\bU(x,t_{n+1}-)$ of \eqref{Euler-1D}-\eqref{data} to the space of piecewise polynomials $\bU(x,t_{n+1}) =\bp_j^{n+1}(x),\ x\in I_j$.

\end{enumerate}

{\bf The Godunov type schemes with the stiffened gas approximation} focuses on the computation of flux functions $\bF_{j+\frac 12}^*$.
Locally, the generalized Riemann problem for \eqref{Euler-1D}  is solved  at each singularity point $(x_{j+\frac 12}, t_n)$, which is shifted to $(0,0)$, subject to the initial data
\begin{equation}
 \ \bU(x,0) =\left\{
 \begin{array}{ll}
\bP_L(x):= \bp_j^{n}(x), \ \ \ &x<0,\\
\bP_R(x):= \bp_{j+1}^{n}(x), & x>0.
 \end{array}
 \right.
 \label{GRP-1}
\end{equation}
The stiffened gas approximation in the last sections can be applied.

Specified to  the second order GRP scheme, the data \eqref{data} takes
\begin{equation}
\bp_j^n(x)= \bU_j^n +(\bU_x)_j^n (x-x_j),  \ \ \ x\in \bI_j,
\end{equation}
where $\bU_j^n$ and $(\bU_x)_j^n$ are constants. Reformulating this data to the form \eqref{data-GRP}, we have
\begin{align}
 &\bU_L =\bU_j^{n}+\frac{\De x}{2}(\bU_x)_j^n,  \ \ \ \ \ \ \bU_L' =(\bU_x)_j^n,\notag\\
 &\bU_R =\bU_{j+1}^{n}-\frac{\De x}{2}(\bU_x)_{j+1}^n, \ \bU_R' =(\bU_x)_{j+1}^n.\notag
 \end{align}
 Then the GRP  numerical fluxes are computed as
\begin{align}
&\bF_{j+\frac 12}^*=\bF\left(\bU^{n+\frac 12}_{j+\frac 12}\right),\ \bU^{n+\frac 12}_{j+\frac 12}:=\bU^{n}_{j+\frac 12}+\frac{\De t}{2}\left(\dfr{\pt \bU}{\pt t}\right)^{n}_{j+\frac 12},\label{2ndflux}
\end{align}
where $\bU^{n}_{j+\frac 12}$ and $\left(\pt \bU/\pt t\right)^{n}_{j+\frac 12}$ are the associated Riemann solution $R^A(0;\bU_L,\bU_R)$ in Section \ref{sec-RP} and $\left(\pt \bU/\pt t \right)_*$ in Section \ref{sec-GRP},
\begin{equation}
\bU^{n}_{j+\frac 12}=R^A(0;\bU_L,\bU_R),  \ \ \ \left(\dfr{\pt \bU}{\pt t}\right)^{n}_{j+\frac 12} =\left(\dfr{\pt \bU}{\pt t}\right)_*.
\end{equation}

\section{Numerical examples}\label{sec:num}

In this section we present several test problems for the Euler equations with general nonlinear equation of states to
display the performance of the  resulting Godunov and GRP schemes. They are three one-dimensional Riemann problems, a two-dimensional Riemann problem, and a shock-bubble problem.
The numerical solutions by the Godunov method with exact EOSs and the Godunov/GRP methods with approximate stiffened gas EOS are denoted in plots as ``Godunov-exact EOS", ``Godunov-approximate EOS" and ``GRP-approximate EOS", respectively. For one dimensional Riemann problems, the exact solution obtained from the exact Riemann solver is shown as a reference solution, denoted as ``exact". The CFL number is chosen as $0.5$ for all cases.

\subsection{1-D examples}

\subsubsection{Contact problem with the Cochran-Chan (C-C) EOS}

\begin{figure}[!htb]
  \includegraphics[width=6.5cm]{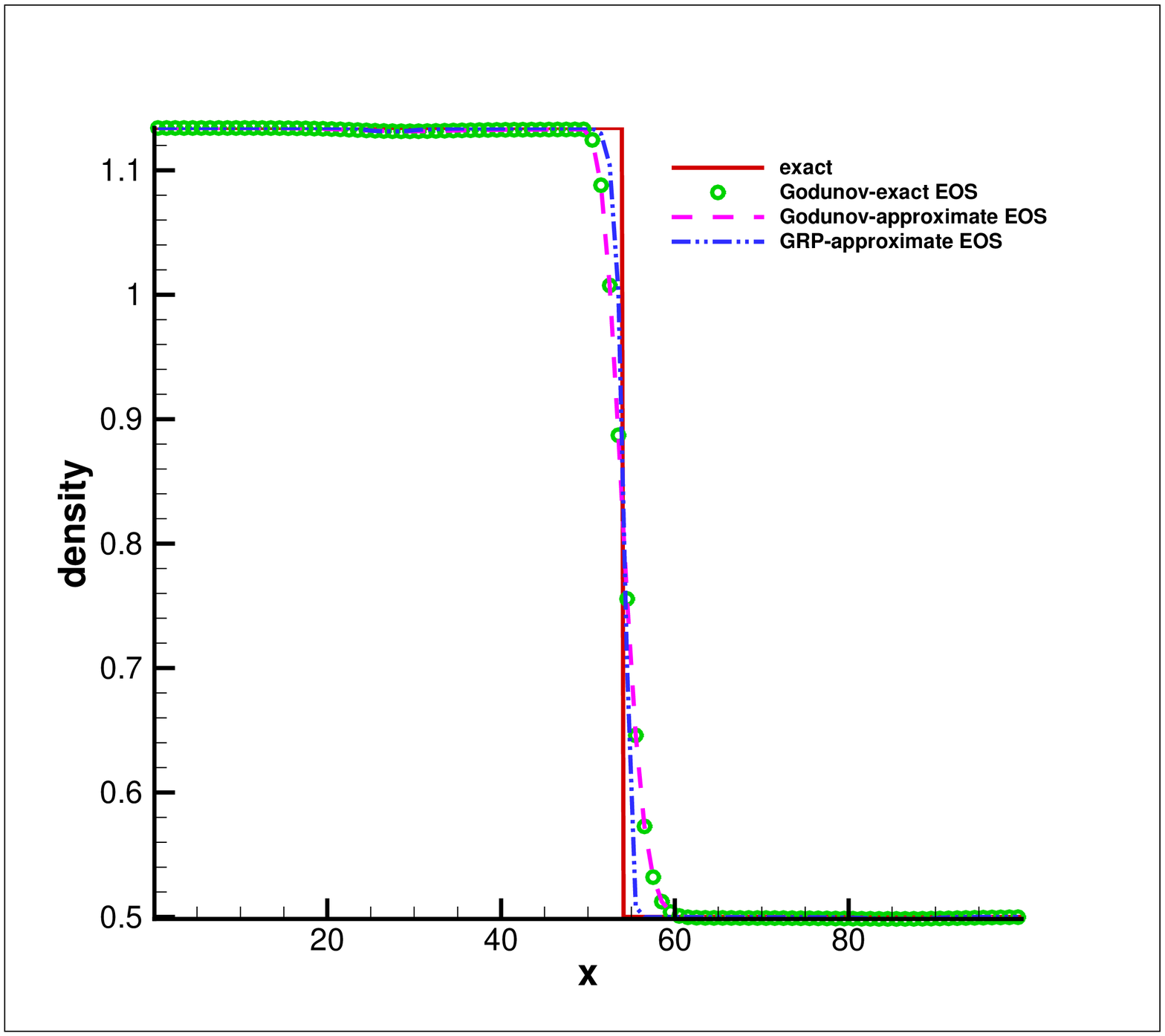}
  \includegraphics[width=6.5cm]{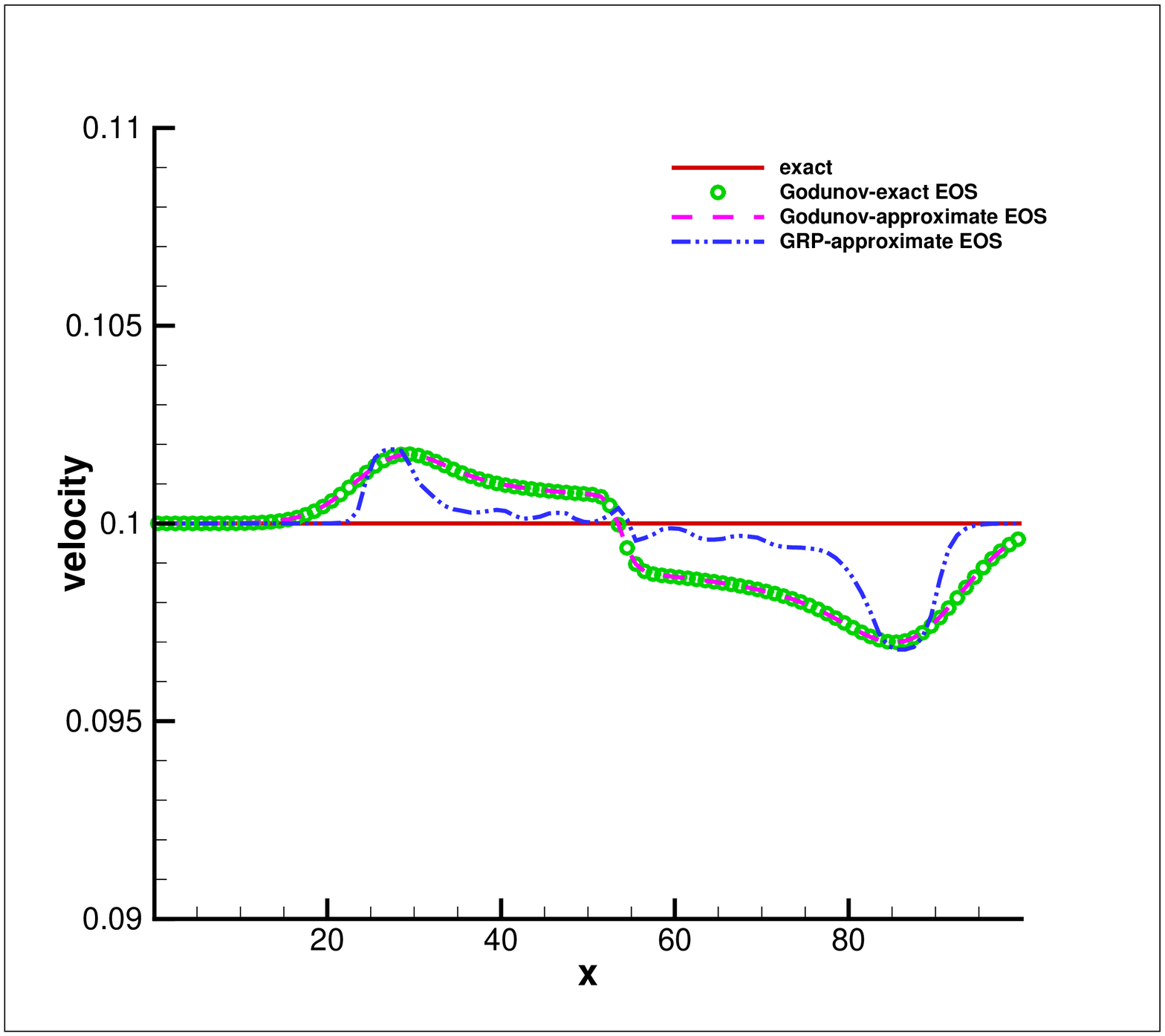}\\
  \includegraphics[width=6.5cm]{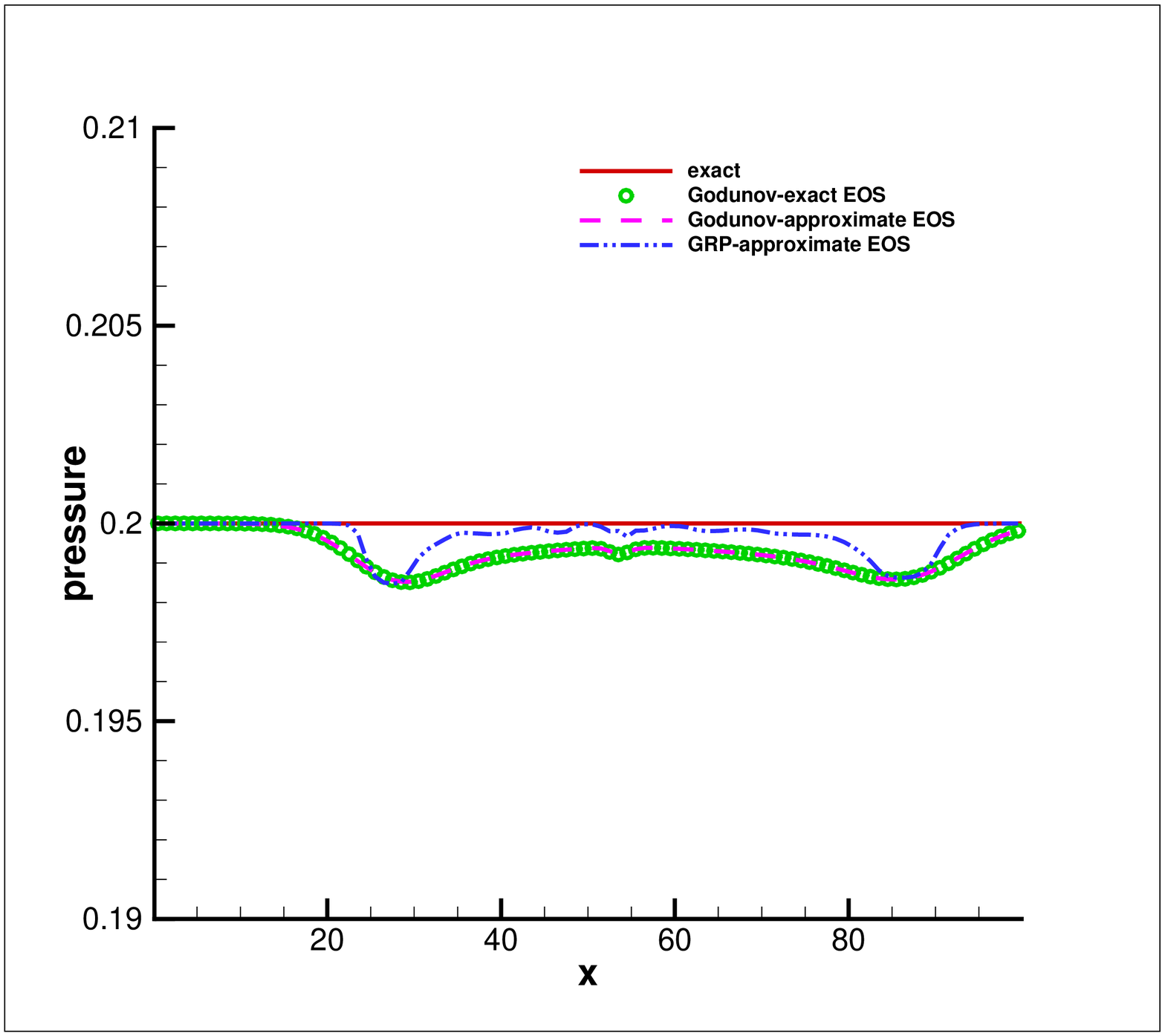}
  \includegraphics[width=6.5cm]{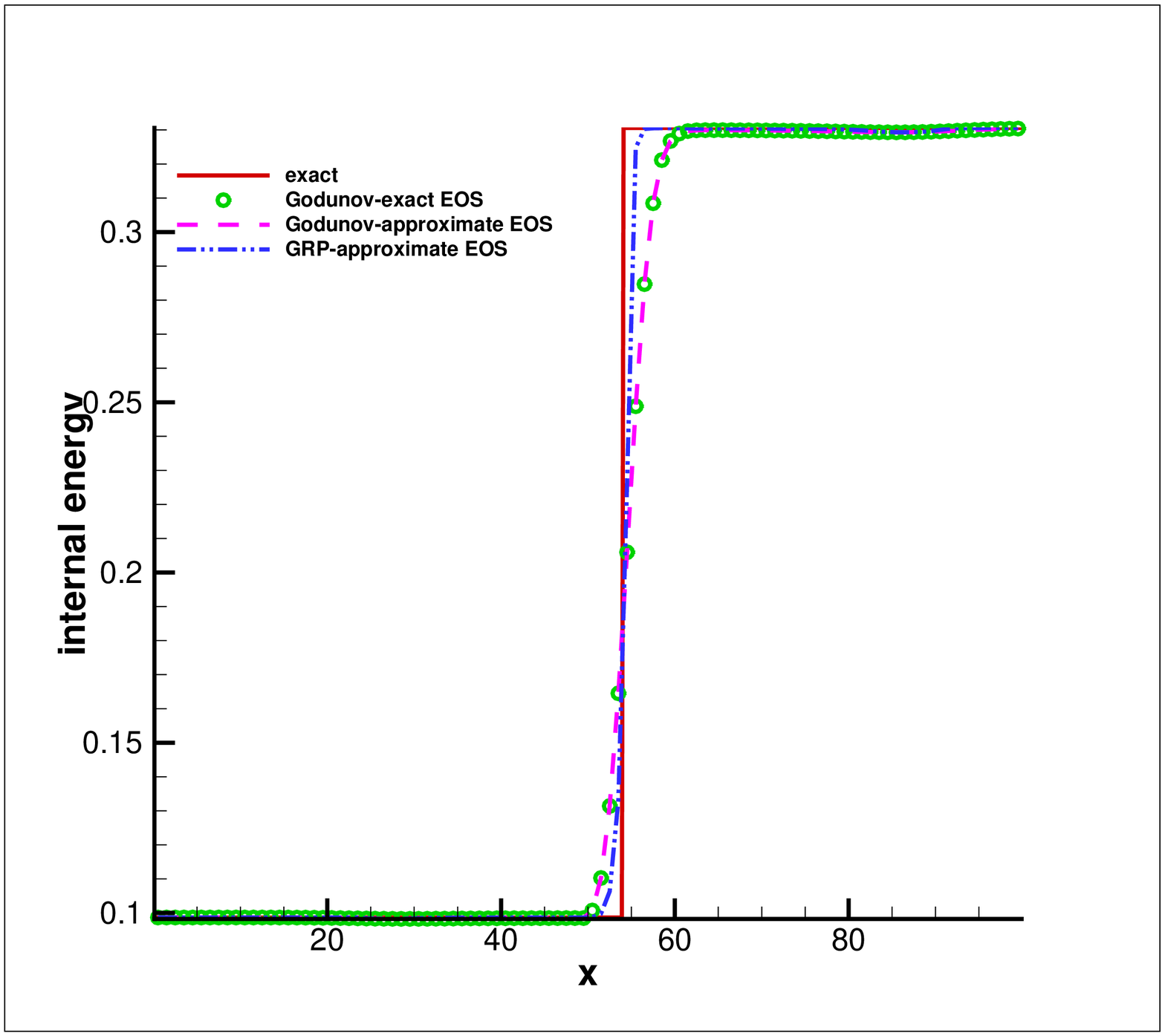}\\
  \caption{Results of density, velocity, pressure and internal energy for the contact problem specified in Table~\ref{contact-param2}. }
  \label{contact-fig1}
\end{figure}

This Riemann problem adopted from \cite{Saurel2007} shows the advection of a density discontinuity in the flow of uniform velocity and
pressure with the C-C EOS,
\begin{align}
&\kappa(\rho)=\Gamma\rho,\ \chi(\rho)=-\Gamma\rho e^{ref}(\rho)+p^{ref}(\rho),\notag\\
&e^{ref}(\rho)=-\frac{A}{(1-\epsilon_1)\rho_0}[(\frac{\rho_0}{\rho})^{1-\epsilon_1}-1]+\frac{B}{(1-\epsilon_2)\rho_0}[(\frac{\rho_0}{\rho})^{1-\epsilon_2}-1]-e_0,\notag\\
&p^{ref}(\rho)=A(\frac{\rho_0}{\rho})^{-\epsilon_1}-B(\frac{\rho_0}{\rho})^{-\epsilon_2}.\notag
\end{align}
 The material parameters for Nitromethane are given in Table~\ref{contact-param1}.
The initial conditions are provided in Table~\ref{contact-param2}.
The computational domain is for $x\in[0,100cm]$.
Figure ~\ref{contact-fig1} shows the numerical solutions and the exact solution for the contact problem using $100$ cells. The numerical solutions are computed using  the Godunov method with exact EOSs and the Godunov/GRP methods with approximate stiffened gas EOS, respectively.
In Figure~\ref{contact-fig1}, we can observe  the contact is precisely captured by the current schemes. The resolutions are almost the same by the exact Riemann solver and the approximate stiffened gas Riemann solver. The approximate GRP method have better resolution thanks to second order accuracy and the inclusion of thermodynamic effect.
We also observe non-physical  pressure and velocity profiles (abeit small deviation) across the contact discontinuity. As discussed in an earlier paper \cite{Lee}, these non-physical phenomena are produced by conservative schemes and arise from  the discrepancy of the physical properties and the nonlinearity of the EOS. However the GRP scheme evidently  improves the result compared with  the first order Godunov scheme.

\begin{table}[!htb]
\centering
\caption{C-C parameters used in the contact problem.}\label{contact-param1}
\begin{tabular}{|c|c|c|c|c|c|c|c|}
  \hline
  C-C & $\rho_0$ & $e_0$ & $\Gamma_0$ & $A$ & $B$ & $\epsilon_1$ & $\epsilon_2$ \\
  EOS & $g/cm^3$ & $Mbar-cm^3/g$ & $-$ & Mbar & Mbar & $-$ & $-$ \\
  \hline
  $\#1$ & 1.134 & 0.0  & 1.19 & 8192  & 1508 & 4.53 &1.42\\
  \hline
\end{tabular}
\end{table}

\begin{table}[!htb]
\centering
\caption{Initial conditions for the contact problem.}\label{contact-param2}
\begin{tabular}{|c|c|c|c|c|c|c|c|c|}
  \hline
  C-C & $x_{intfc}^0$ & $t_{fin}$ & $\rho_L$ & $p_L$ & $u_L$ & $\rho_R$ & $p_R$ & $u_R$ \\
  EOS &  $cm$ & $\mu s$ & $g/cm^3$ & $Mbar$ & $cm/\mu s$ & $g/cm^3$ & $Mbar$ & $cm/\mu s$ \\
  \hline
  $\#1$& 50 & 40 &1.134& 2e4& 0.1 &0.5 &2e4 &0.1\\
  \hline
\end{tabular}
\end{table}

\subsubsection{Shyue Shock Tube}
This problem, proposed by Shyue \cite{Shyue} as a test problem for algorithm development, uses
the JWL EOS
  \begin{align}
  &\kappa(\rho)=\Gamma\rho,\ \chi(\rho)=-\Gamma\rho e^{ref}(\rho)+p^{ref}(\rho),\label{Mie}\\
  &e^{ref}(\rho)=\frac{A}{R_1\rho_0}\exp(-R_1\frac{\rho_0}{\rho})+\frac{B}{R_2\rho_0}\exp(-R_2\frac{\rho_0}{\rho})-e_0,\label{JWL-1}\\
  &p^{ref}(\rho)=A\exp(-R_1\frac{\rho_0}{\rho})+B\exp(-R_2\frac{\rho_0}{\rho}),\label{JWL-2}
  \end{align}
on either side of a massless discontinuity. The JWL parameters,
listed in Table~\ref{shyue-param1}, are an approximation for the EOS of TNT. This problem is set over the
interval $0\leqslant x\leqslant100 cm$, with the initial conditions  in Table~\ref{shyue-param2}.

\begin{figure}[!htb]
  \includegraphics[width=6.5cm]{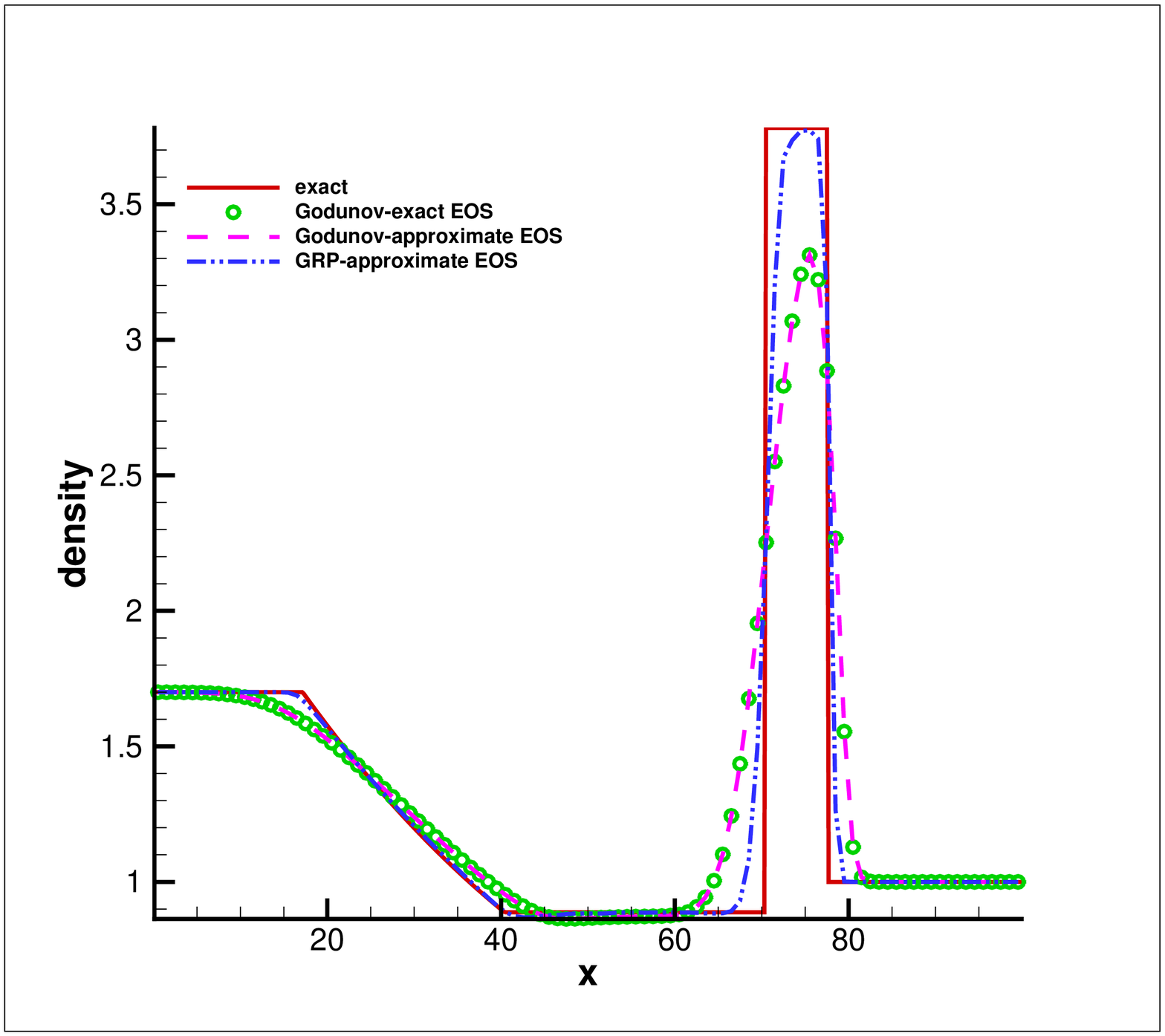}
  \includegraphics[width=6.5cm]{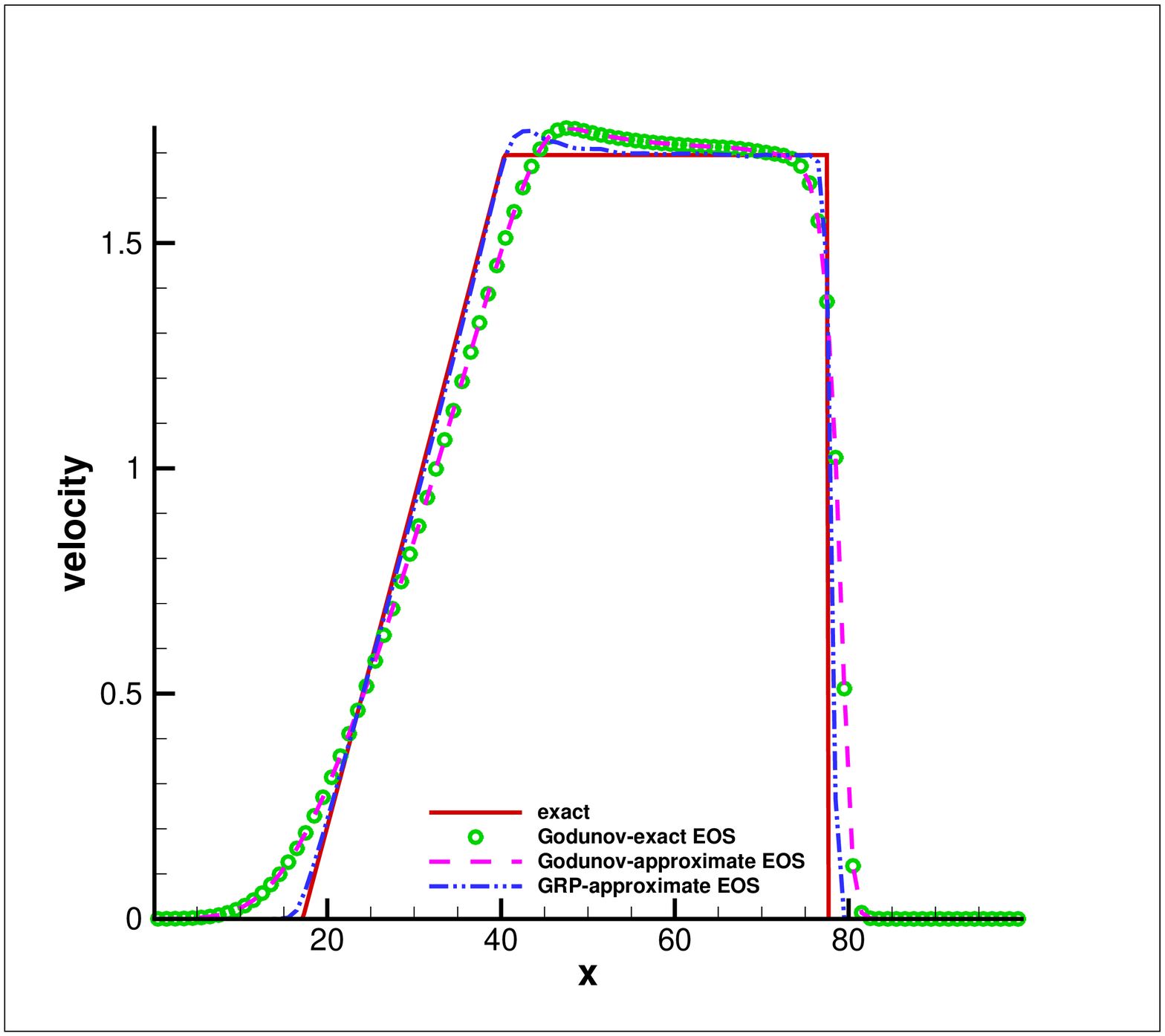}\\
  \includegraphics[width=6.5cm]{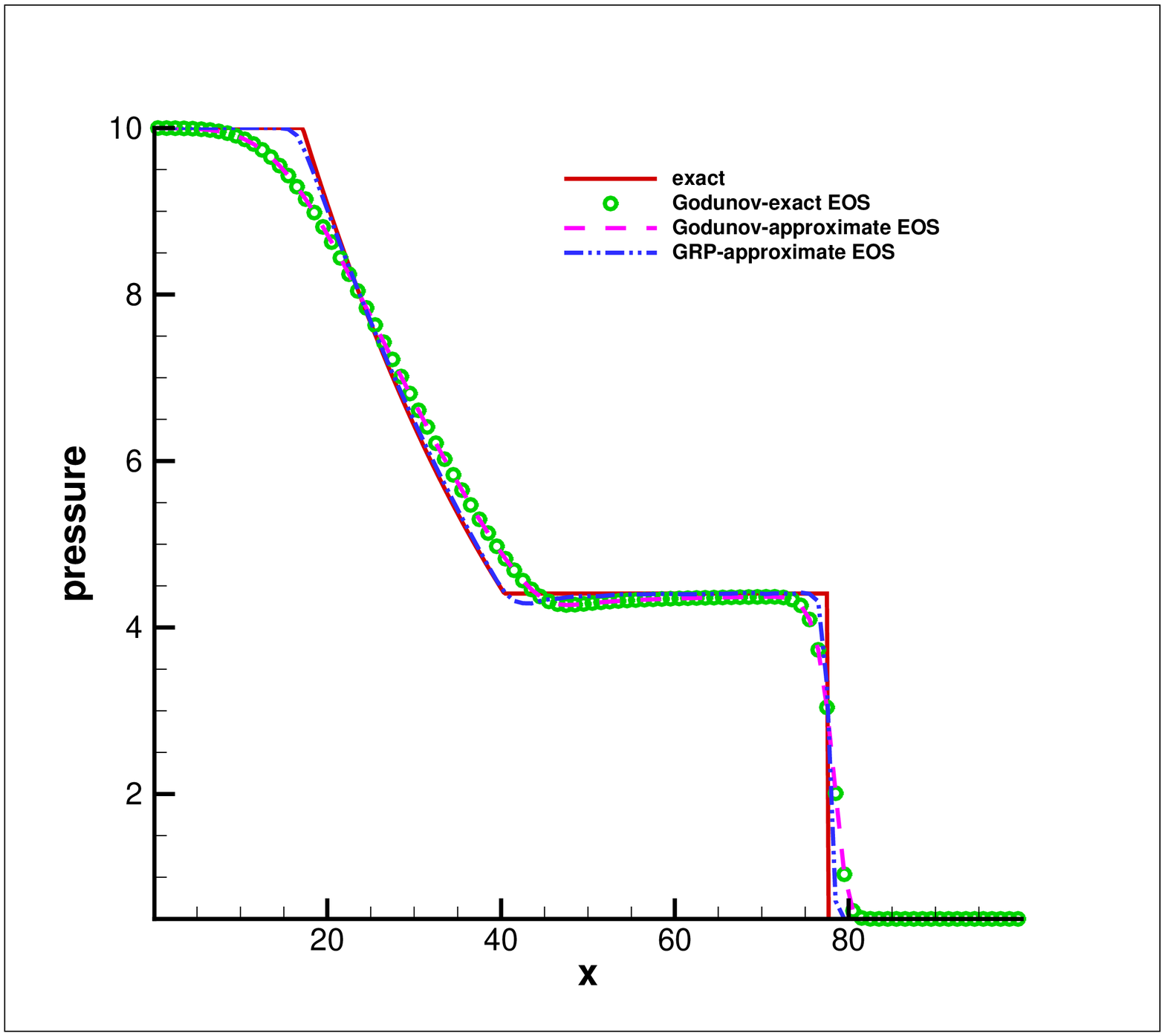}
  \includegraphics[width=6.5cm]{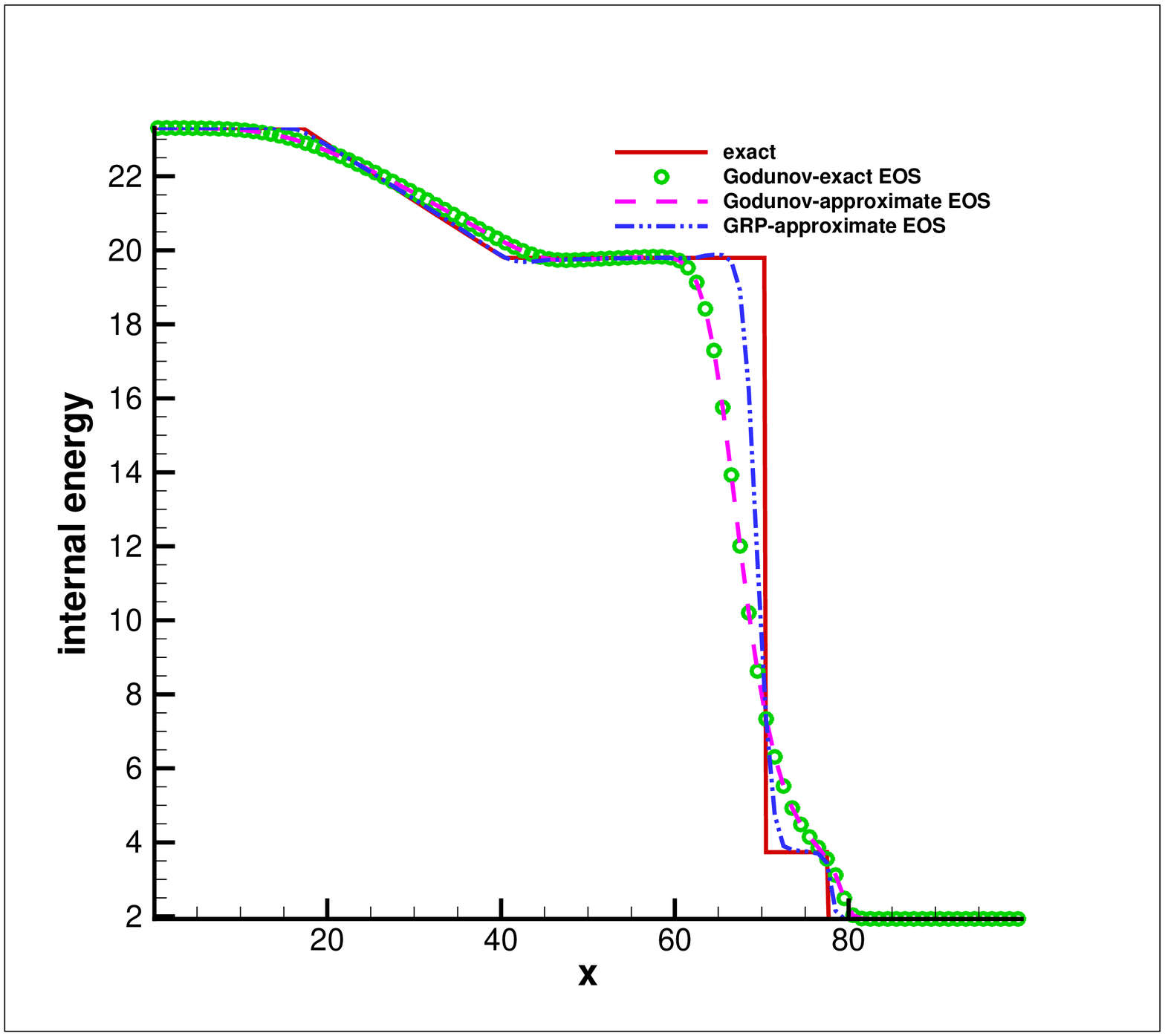}\\
  \caption{Results of density, velocity, pressure and internal energy for the Shyue shock tube problem specified in Table~\ref{shyue-param2}. }
  \label{shyue-fig2}
\end{figure}

Results for this problem at $t=12\mu s$  are given in Figure ~\ref{shyue-fig2} with $100$ computational cells, which contains snapshot plots of the
density, velocity, pressure, and internal energy. The red line represents results obtained with  the general
EOS exact solution code, the green circles are obtained by the Godunov scheme with the exact Riemann solver, and the purple and blue lines are obtained  by the approximate Godunov and GRP schemes using the proposed stiffened gas approximation, respectively.
The numerical results computed by the Godunov scheme with approximate stiffened gas EOS show the same resolution compared with the one with exact JWL EOS.
The approximate GRP scheme shows better resolution with second order accuracy.

\begin{table}[!htb]
\centering
\caption{JWL parameters for TNT used in the Shyue shock tube problem.}\label{shyue-param1}
\begin{tabular}{|c|c|c|c|c|c|c|c|c|}
  \hline
  JWL & $\rho_0$ & $e_0$ & $\Gamma_0$ & $A$ & $B$ & $R_1$ & $R_2$ \\
  EOS & $g/cm^3$ & $Mbar-cm^3/g$ & $-$ & Mbar & Mbar & $-$ & $-$ \\
  \hline
  $\#2$ & 1.84 & 0.0  & 0.25 & 8.545  &0.205& 4.6 &1.35\\
  \hline
\end{tabular}
\end{table}

\begin{table}[!htb]
\centering
\caption{Initial conditions for the Shyue shock tube problem.}\label{shyue-param2}
\begin{tabular}{|c|c|c|c|c|c|c|c|c|}
  \hline
  JWL & $x_{intfc}^0$ & $t_{fin}$ & $\rho_L$ & $p_L$ & $u_L$ & $\rho_R$ & $p_R$ & $u_R$ \\
  EOS &  $cm$ & $\mu s$ & $g/cm^3$ & $Mbar$ & $cm/\mu s$ & $g/cm^3$ & $Mbar$ & $cm/\mu s$ \\
  \hline
  $\#2$& 50 & 12 &1.7& 10.0& 0.0 &1.0 &0.5 &0.0\\
  \hline
\end{tabular}
\end{table}

%


\subsubsection{Lee Shock Tube}
Lee et al. \cite{Lee} consider another JWL 1D Riemann problem. As in the previous problem, the same JWL EOS is used on
either side of the initial discontinuity. The JWL parameters, listed in Table~\ref{lee-param1}, are an
approximation for the EOS of LX-17. Like the previous example, this problem is also
run on the interval $0\leqslant x\leqslant100 cm$, but with the different initial conditions, as given
in Table~\ref{lee-param2}.

\begin{figure}[!htb]
  \includegraphics[width=6.5cm]{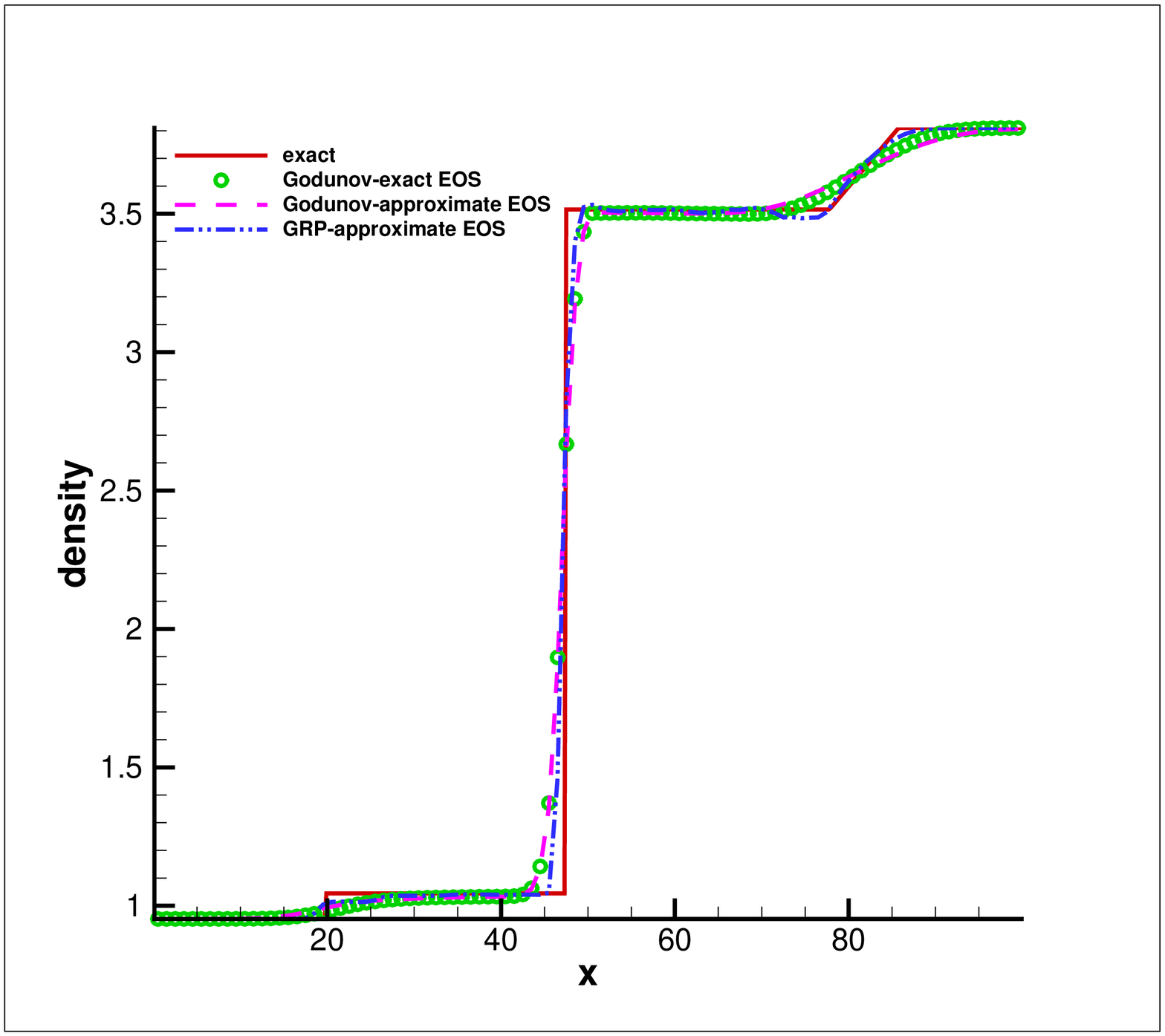}
  \includegraphics[width=6.5cm]{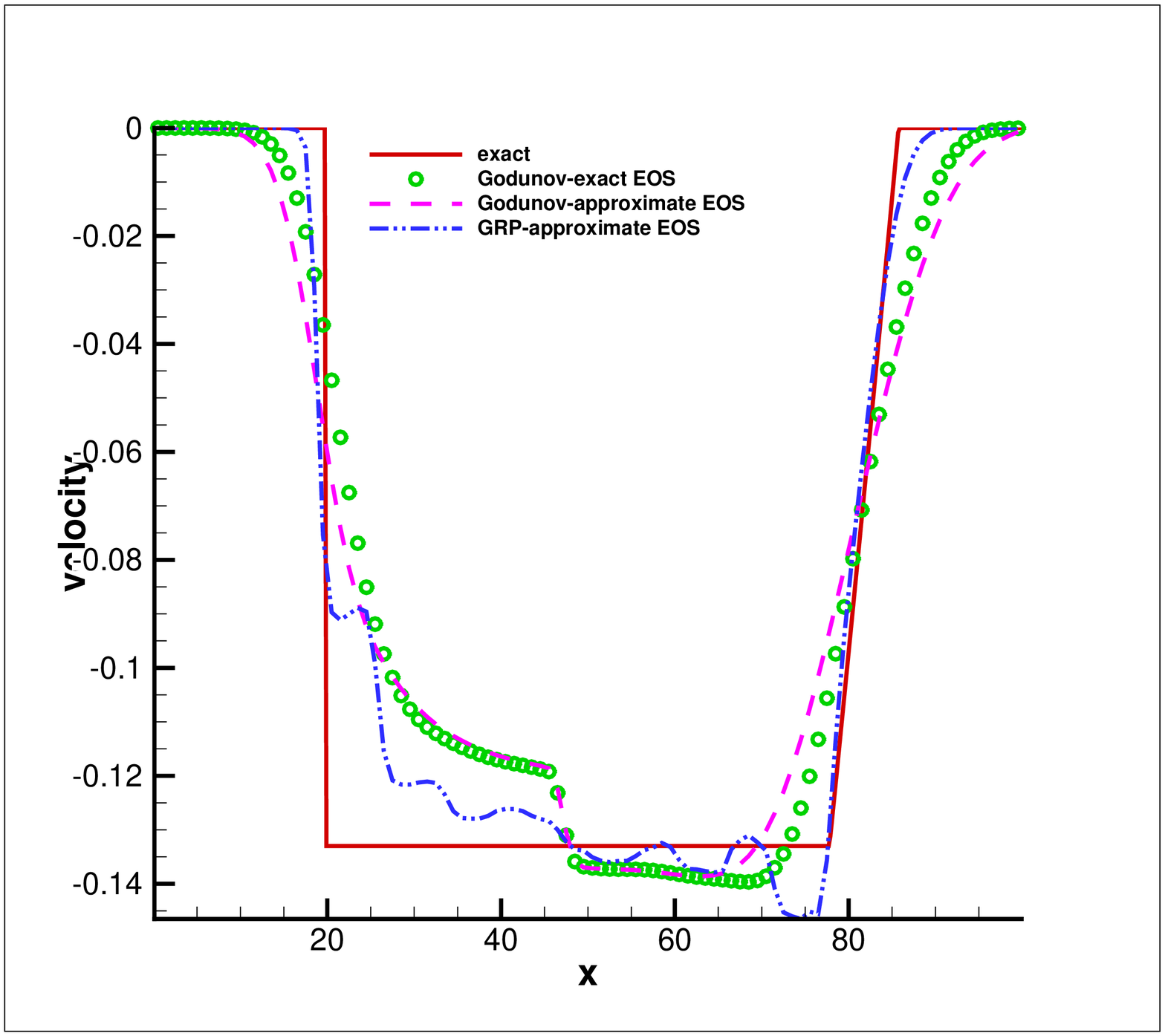}\\
  \includegraphics[width=6.5cm]{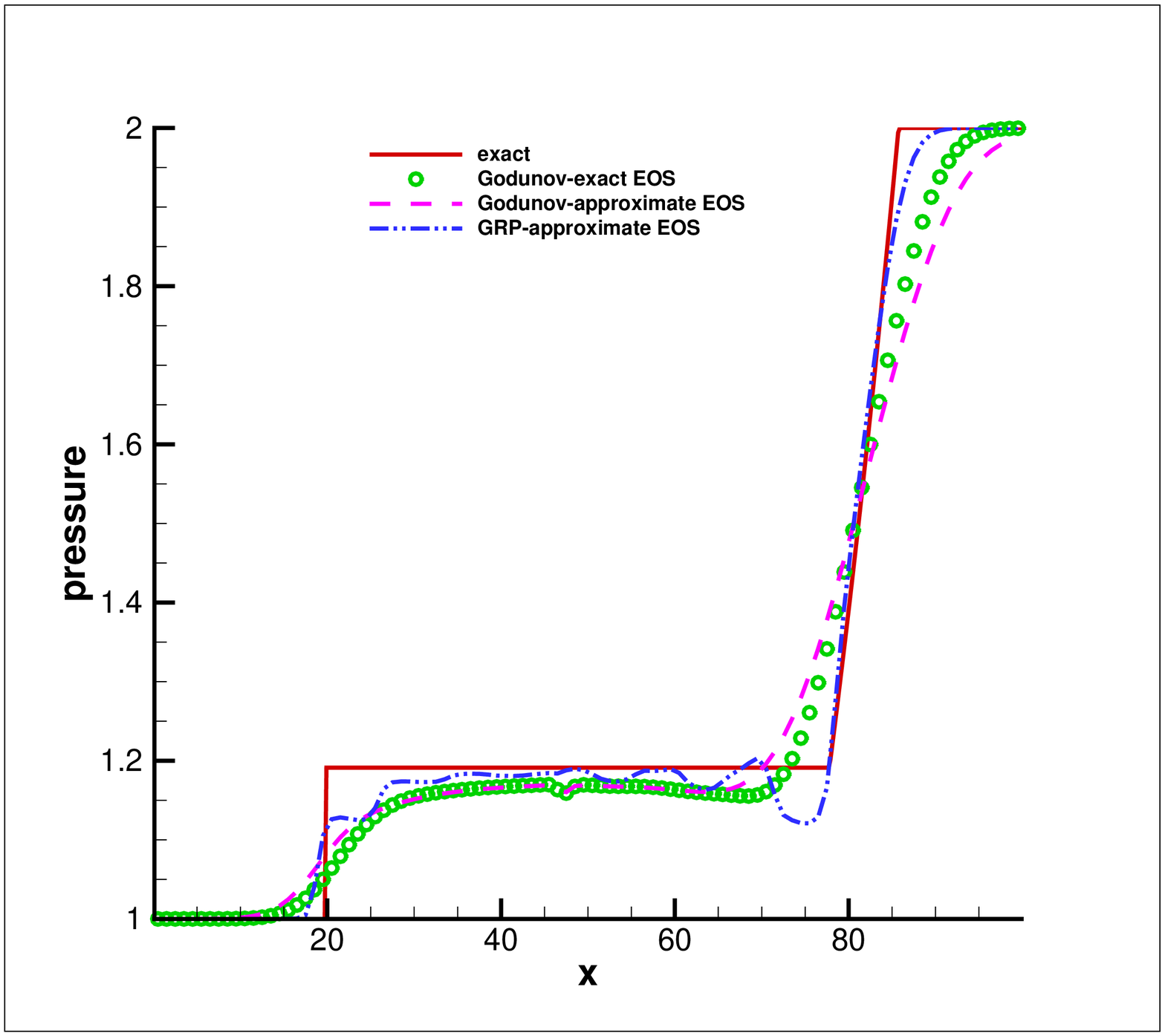}
    \includegraphics[width=6.5cm]{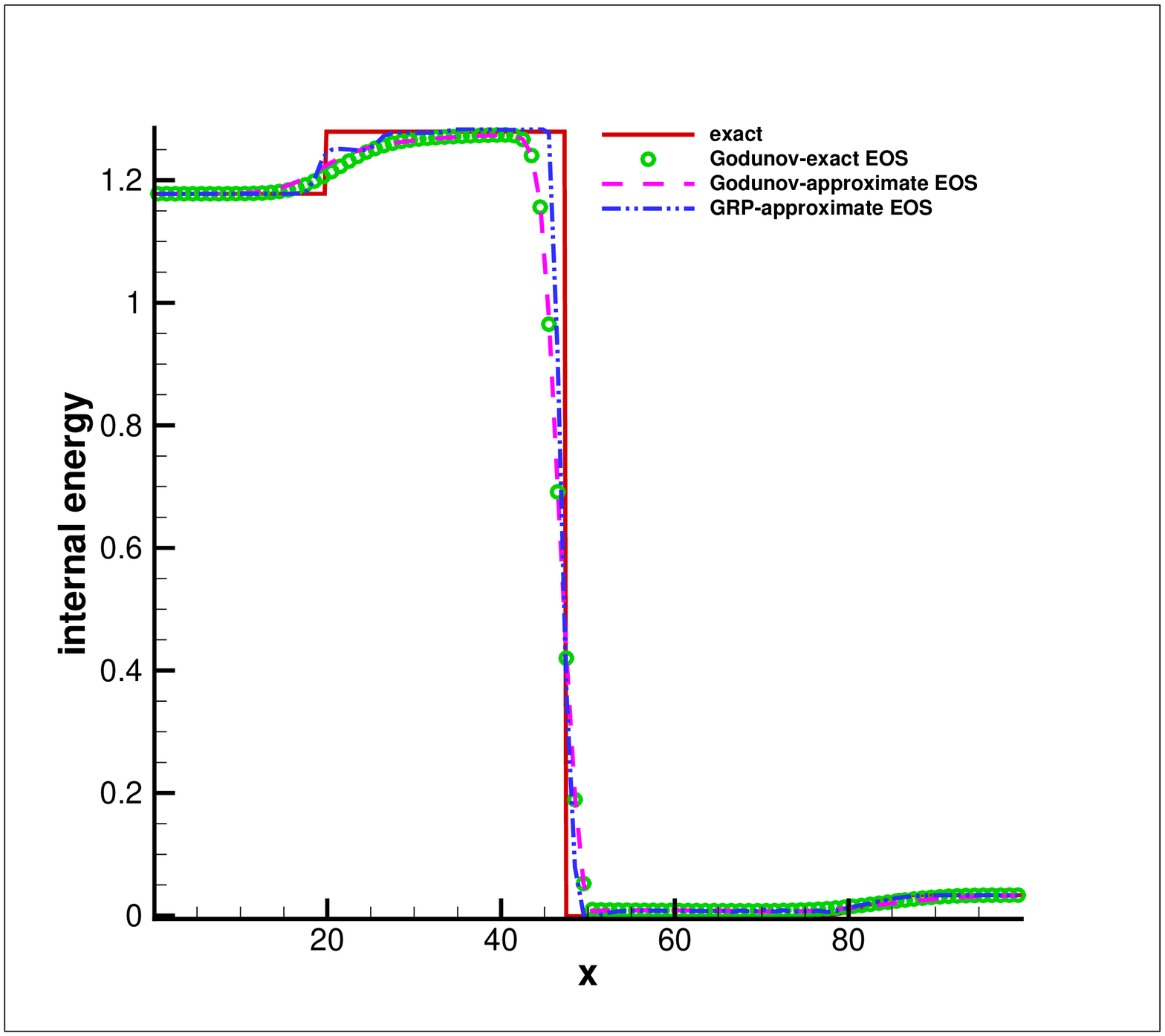}\\
  \caption{Results of density, velocity, pressure and internal energy for the Lee shock tube problem specified in Table~\ref{lee-param2}. }\label{lee-fig2}
\end{figure}

Results for this problem are given in Figure ~\ref{lee-fig2}, which contains snapshot plots of the
density, velocity, pressure, and internal energy.
Numerical results computed by the Godunov and GRP schemes with approximate stiffened gas EOS show reasonable resolution.
Spurious oscillation in pressure and velocity appears due to the sufficiently nonlinear
EOS \cite{banks,Lee,Saurel2007}.

\begin{table}[!htb]
\centering
\caption{JWL parameters  in the Lee shock tube problem.}\label{lee-param1}
\begin{tabular}{|c|c|c|c|c|c|c|c|c|}
  \hline
  JWL & $\rho_0$ & $e_0$ & $\Gamma_0$ & $A$ & $B$ & $R_1$ & $R_2$ \\
  EOS & $g/cm^3$ & $Mbar-cm^3/g$ & $-$ & Mbar & Mbar & $-$ & $-$ \\
  \hline
  $\#3$ & 1.905& 0.0 &0.8938& 632.1& -0.04472& 11.3& 1.13\\
  \hline
\end{tabular}
\end{table}

\begin{table}[!htb]
\centering
\caption{Initial conditions for the Lee shock tube problem.}\label{lee-param2}
\begin{tabular}{|c|c|c|c|c|c|c|c|c|}
  \hline
  JWL & $x_{intfc}^0$ & $t_{fin}$ & $\rho_L$ & $p_L$ & $u_L$ & $\rho_R$ & $p_R$ & $u_R$ \\
  EOS &  $cm$ & $\mu s$ & $cm^3/g$ & $Mbar$ & $cm/\mu s$ & $cm^3/g$ & $Mbar$ & $cm/\mu s$ \\
  \hline
  $\#3$& 50 & 20 & 0.9525 &1.0& 0.0 &3.810& 2.0& 0.0\\
  \hline
\end{tabular}
\end{table}

\subsection{2-D examples}

\subsubsection{2-D Riemann problem}
\begin{figure}[!htb]
  \includegraphics[width=6.5cm]{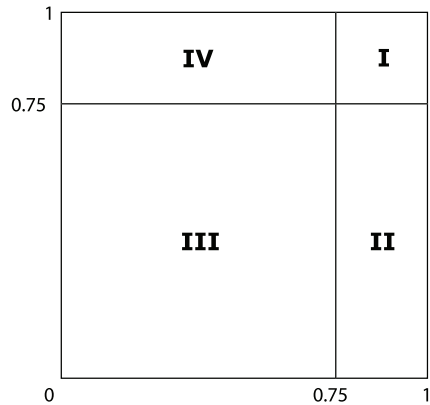}\\
  \caption{Two-dimensional Riemann problem: Initial domain}\label{2DRP-fig2}
\end{figure}

\begin{figure}[!htb]
  \includegraphics[width=6.5cm]{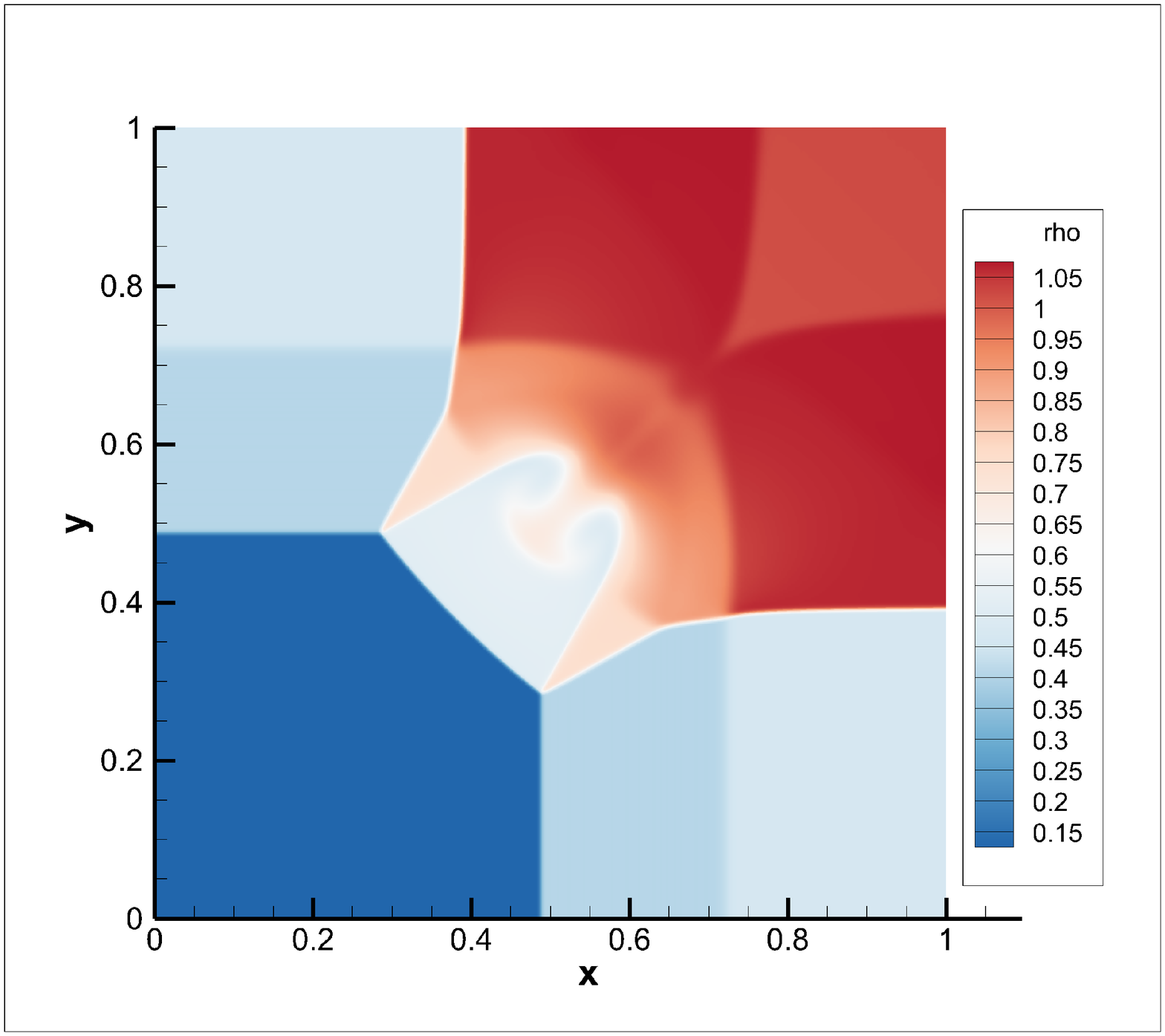}
  \includegraphics[width=6.5cm]{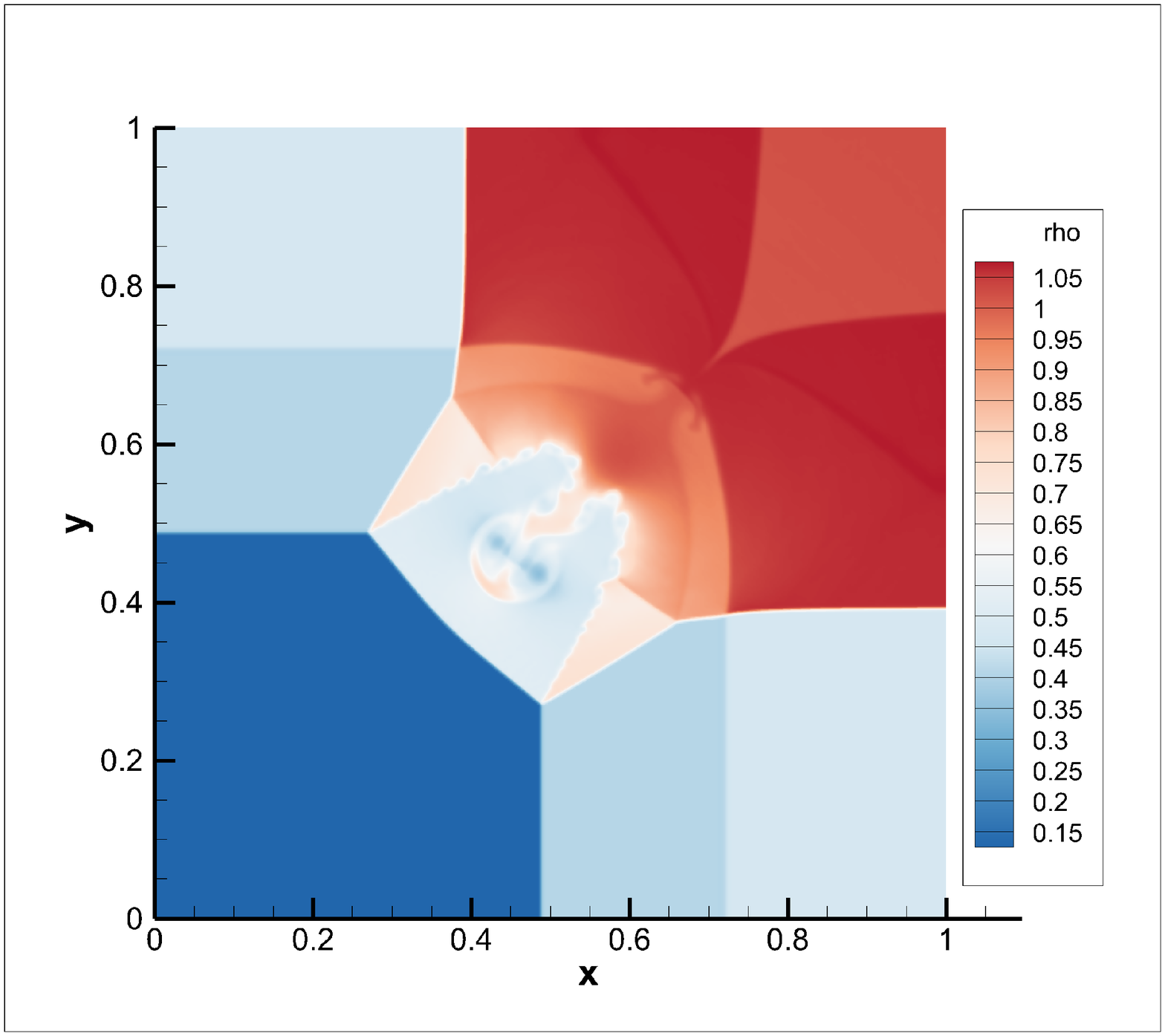}\\
  \includegraphics[width=6.5cm]{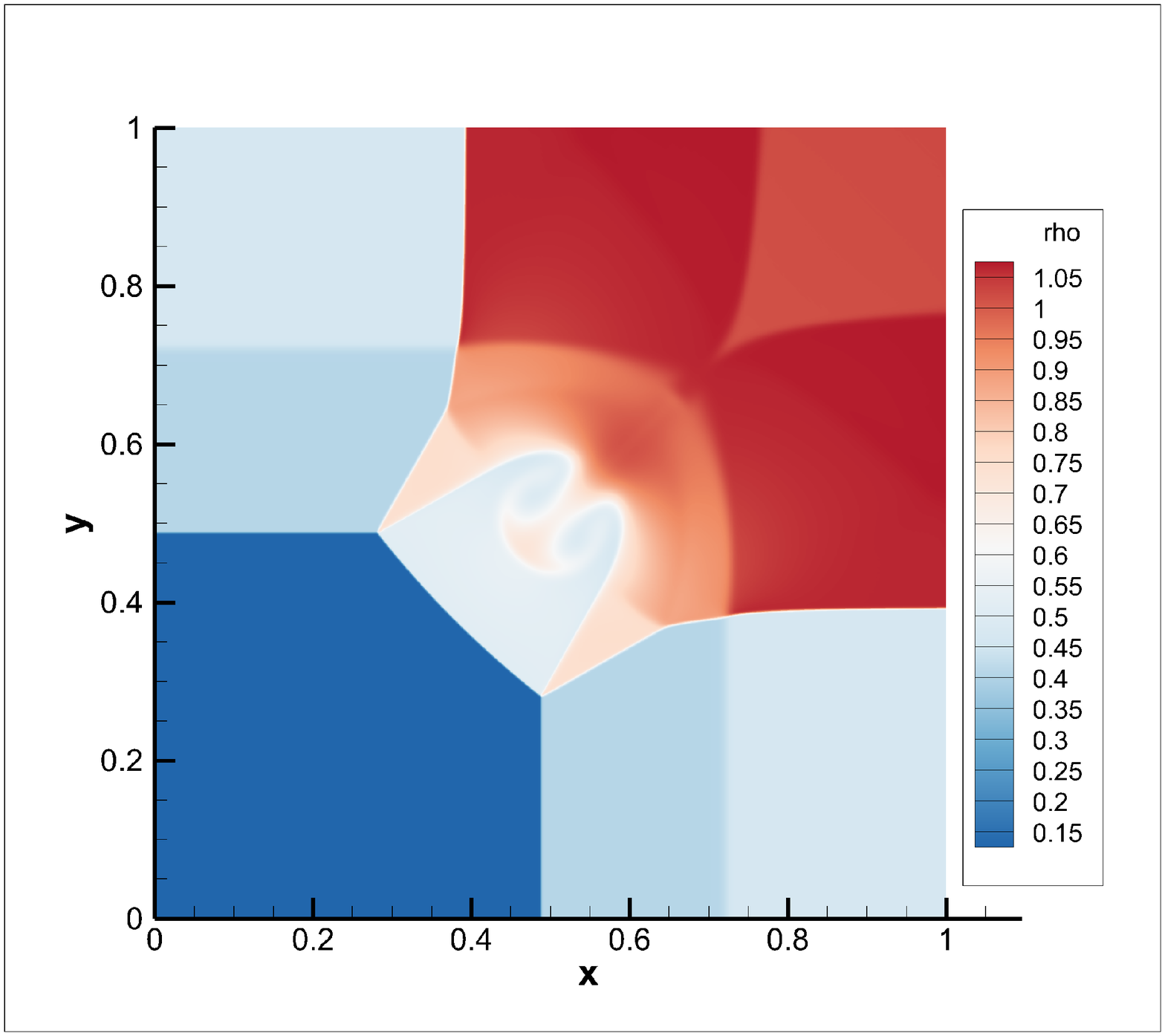}
  \includegraphics[width=6.5cm]{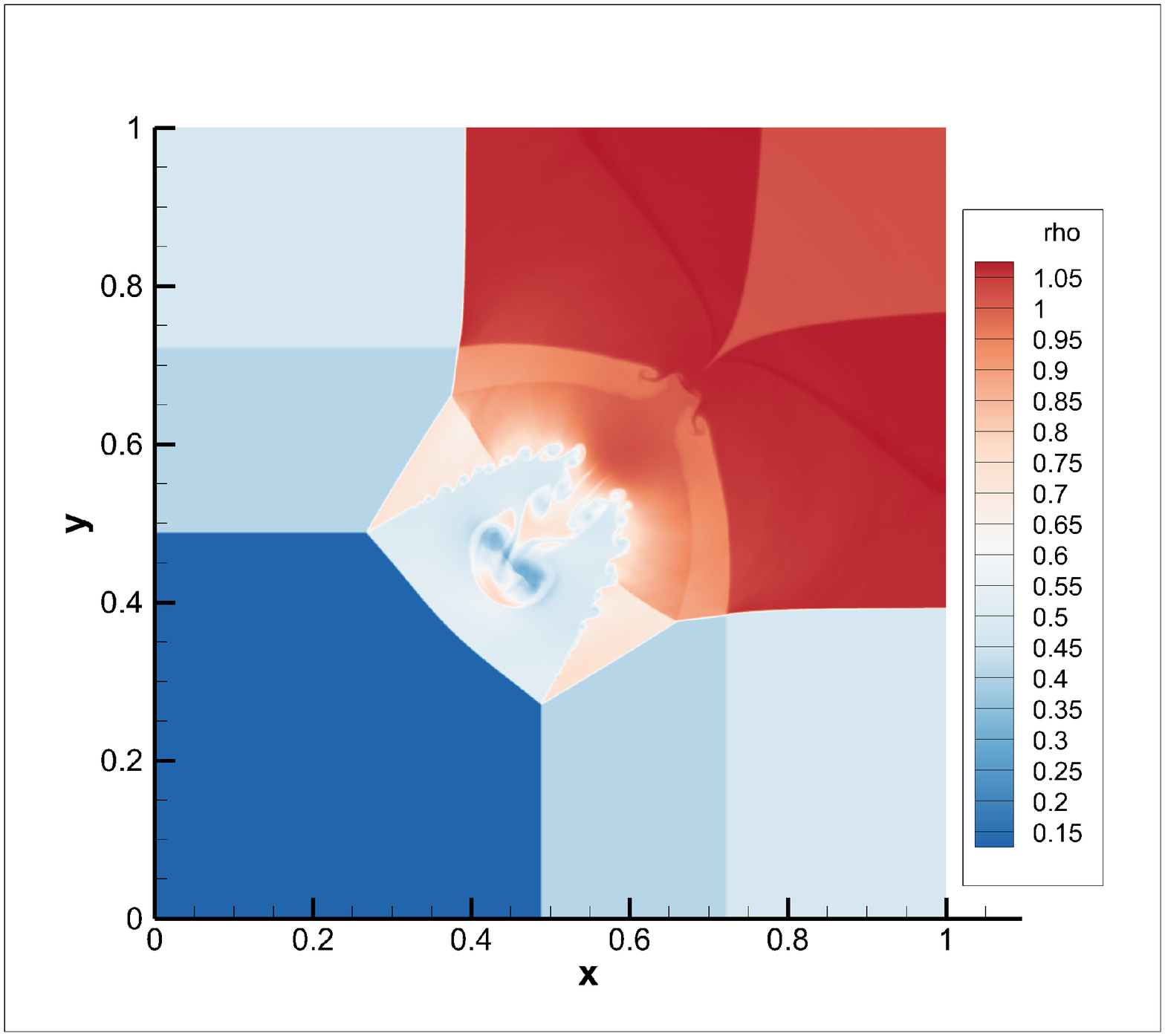}\\
  \caption{ The density contours for two-dimensional Riemann problem: computed by the Godunov (left) and GRP (right) methods with the stiffened gas approximation (Upper: $400\times 400$ cells; Lower:$800 \times 800$ cells)}\label{2DRP-fig1}
\end{figure}

The two-dimensional test case consists of four initial states divided by the horizontal and vertical diaphragms, as
shown in Figure~ \ref{2DRP-fig2} \cite{Bai}.
Here, just for convenience, we use non-dimensional variables
denoted with $\tilde{()}$. Variables are normalised by a reference pressure of $100 GPa$, a reference density chosen to be $\rho_0$ in Table ~\ref{2drp-param1}, a reference length of $0.1 m$, and a reference time of $10^{-4} s$.
For the domain of $1.0\times1.0$ in $\tilde{x}-\tilde{y}$ space, the intersection of two diaphragms is located at $(0.75,0.75)$.
The initial conditions for each region are summarised in Table.~\ref{2drp-param1}. The JWL EOS for the detonation product of PBX-9502 is used. The numerical solutions of the Riemann problem at $\tilde{t} = 0.5$ using the Godunov and GRP schemes with approximated stiffened gas EOS are shown in Figure~ \ref{2DRP-fig1} in the density contours with equidistant computational cells of $400\times 400$ (upper) and $800\times800$ (lower). We can observe their respective performances clearly.


\begin{table}[!htb]
\centering
\caption{JWL parameters for PBX-9502 used in the two dimensional Riemann problem.}\label{2drp-param1}
\begin{tabular}{|c|c|c|c|c|c|c|c|c|}
  \hline
  JWL & HE & $\rho_0$ & $e_0$ & $\Gamma_0$ & $A$ & $B$ & $R_1$ & $R_2$ \\
  EOS &   & $g/cm^3$ & $Mbar-cm^3/g$ & $-$ & Mbar & Mbar & $-$ & $-$ \\
  \hline
  $\#4$& PBX-9502 & 1.895 & 0 & 0.5 & 1.362e6& 7.199e4& 6.2 &2.2\\
  \hline
\end{tabular}
\end{table}

\begin{table}[!htb]
\centering
\caption{Initial condition for the two-dimensional Riemann problem.}\label{2drp-param2}
\begin{tabular}{|c|c|c|c|c|}
  \hline
  JWL EOS $\#4$ & Region \uppercase\expandafter{\romannumeral1} & Region \uppercase\expandafter{\romannumeral2} & Region \uppercase\expandafter{\romannumeral3} & Region \uppercase\expandafter{\romannumeral4} \\
  \hline
  $\tilde{\rho}$& 1 & 0.5 & 0.125 &0.5\\
  \hline
  $\tilde{u}$& 0 & 0 & 1 &1\\
  \hline
  $\tilde{v}$& 0 & 1 & 1 &0\\
  \hline
  $\tilde{p}$& 1 & 0.25 & 0.025 &0.25\\
  \hline
\end{tabular}
\end{table}

\subsubsection{2-D shock-bubble interaction with JWL EOS}

This last numerical experiment shows the interaction between a right shock and a bubble.
The initial condition used for the shock-bubble interaction case is summarised in Table ~\ref{bubble-param1}.
For the domain of $[0,300]cm\times[0,100]cm$, the shocked state is set for $x<5cm$. The radius of the bubble is $r=25cm$ and the center of it is located at $(50cm, 50cm)$.
The JWL EOS with parameters for TNT is applied for this case (Table ~\ref{shyue-param1}).
\begin{table}[!htb]
\centering
\caption{Initial conditions for the shock-bubble interaction problem.}\label{bubble-param1}
\begin{tabular}{|c|c|c|c|}
  \hline
  states & shocked & unshocked & bubble \\
 \hline
  $\rho(g/cm^3)$ & 4.545 & 1 & 2\\
  \hline
  $u(cm/\mu s)$ & 1.737 & 0 & 0\\
  \hline
  $v(cm/\mu s)$ & 0 & 0 & 0\\
  \hline
  $p(Mbar)$ & 4.369 & 0.5 & 0.5\\
  \hline
\end{tabular}
\end{table}
The numerical solutions of the shock-bubble interaction problem at $t = 40\mu s$ and $t=70\mu s$ using the Godunov scheme  the GRP scheme, respectively, with approximated stiffened gas EOS are shown in Figure~ \ref{bubble-fig1} in the density contours equidistant computational cells of $1200\times 400$ (upper) and $2400\times800$ (lower).
It can be seen that the second order GRP scheme with approximate EOS shows better resolution and small structures compared with the Godunov scheme.


\begin{figure}[!htb]
  \includegraphics[width=7cm]{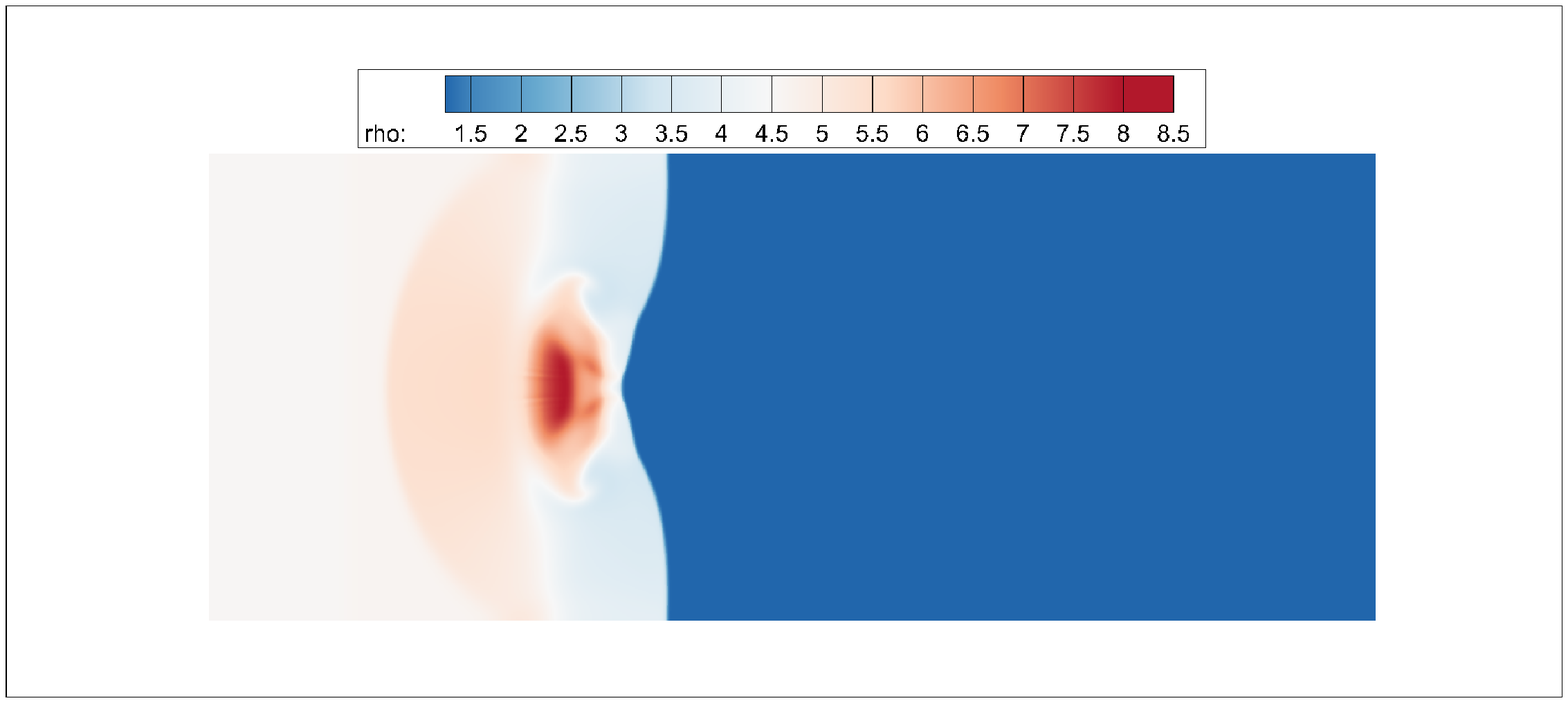}
  \includegraphics[width=7cm]{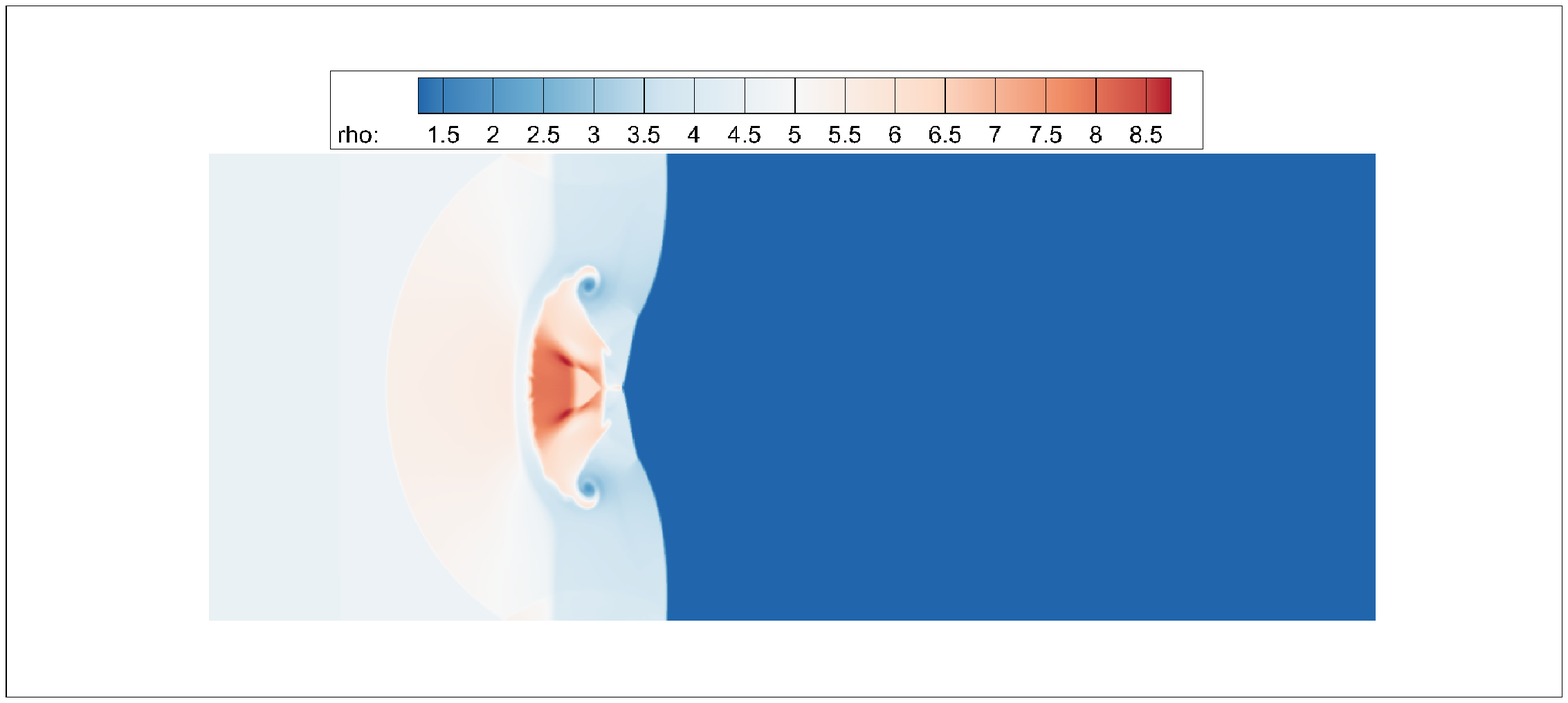}\\
  \includegraphics[width=7cm]{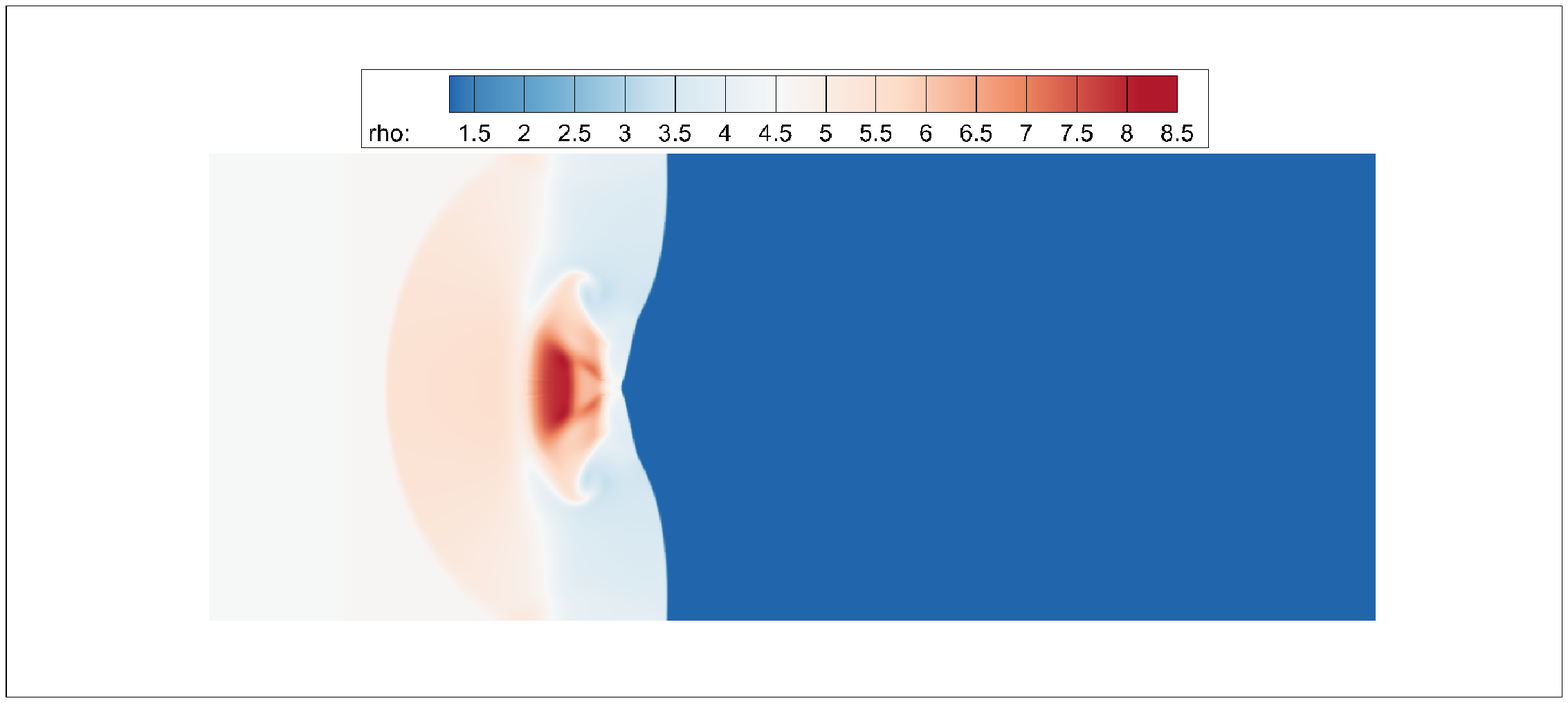}
    \includegraphics[width=7cm]{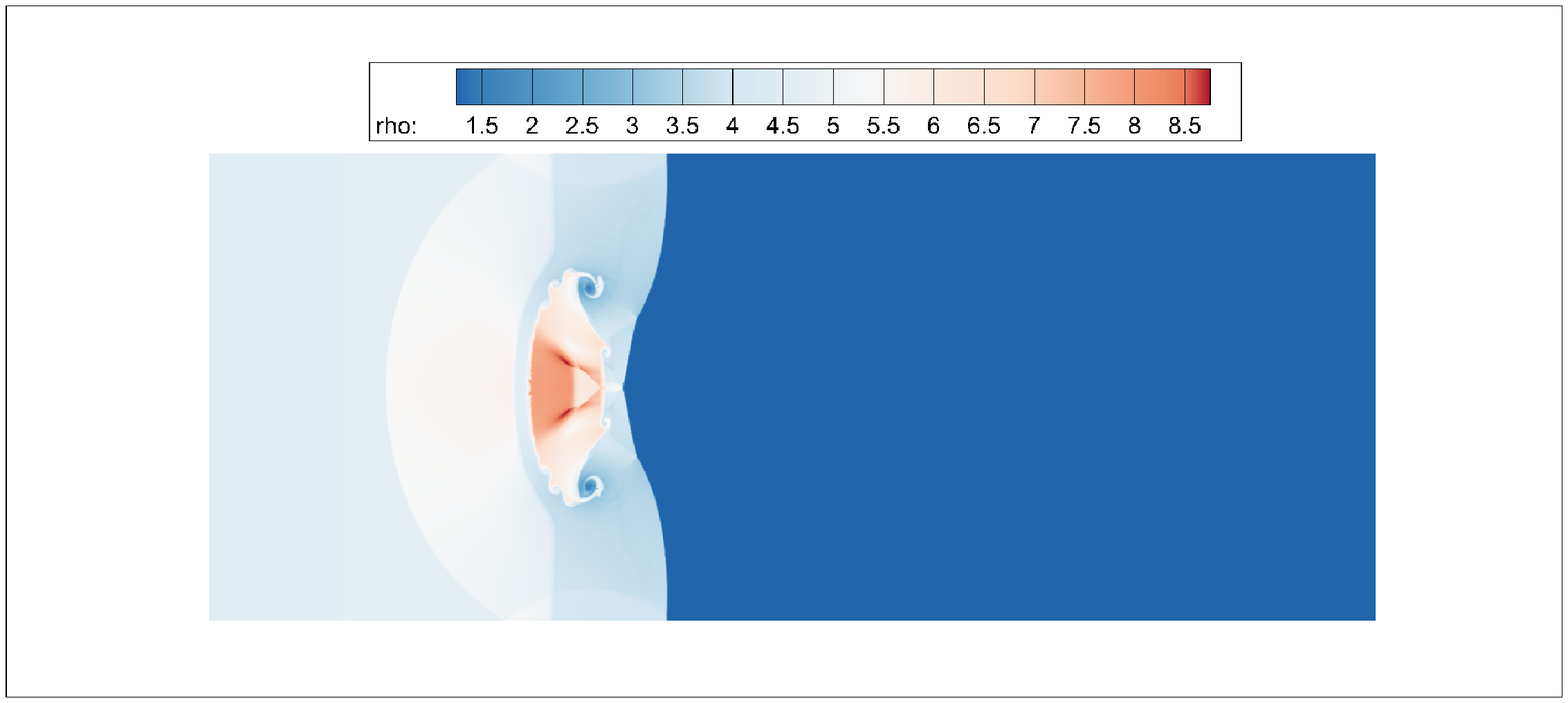}\\
  \caption{The density contours for shock-bubble interaction problem at $t = 40\mu s$: computed by the Godunov (left) and GRP (right) methods with the stiffened gas approximation (Upper: $1200\times 400$ cells; Lower:$2400 \times 800$ cells). }\label{bubble-fig1}
\end{figure}

\begin{figure}[!htb]
  \includegraphics[width=7cm]{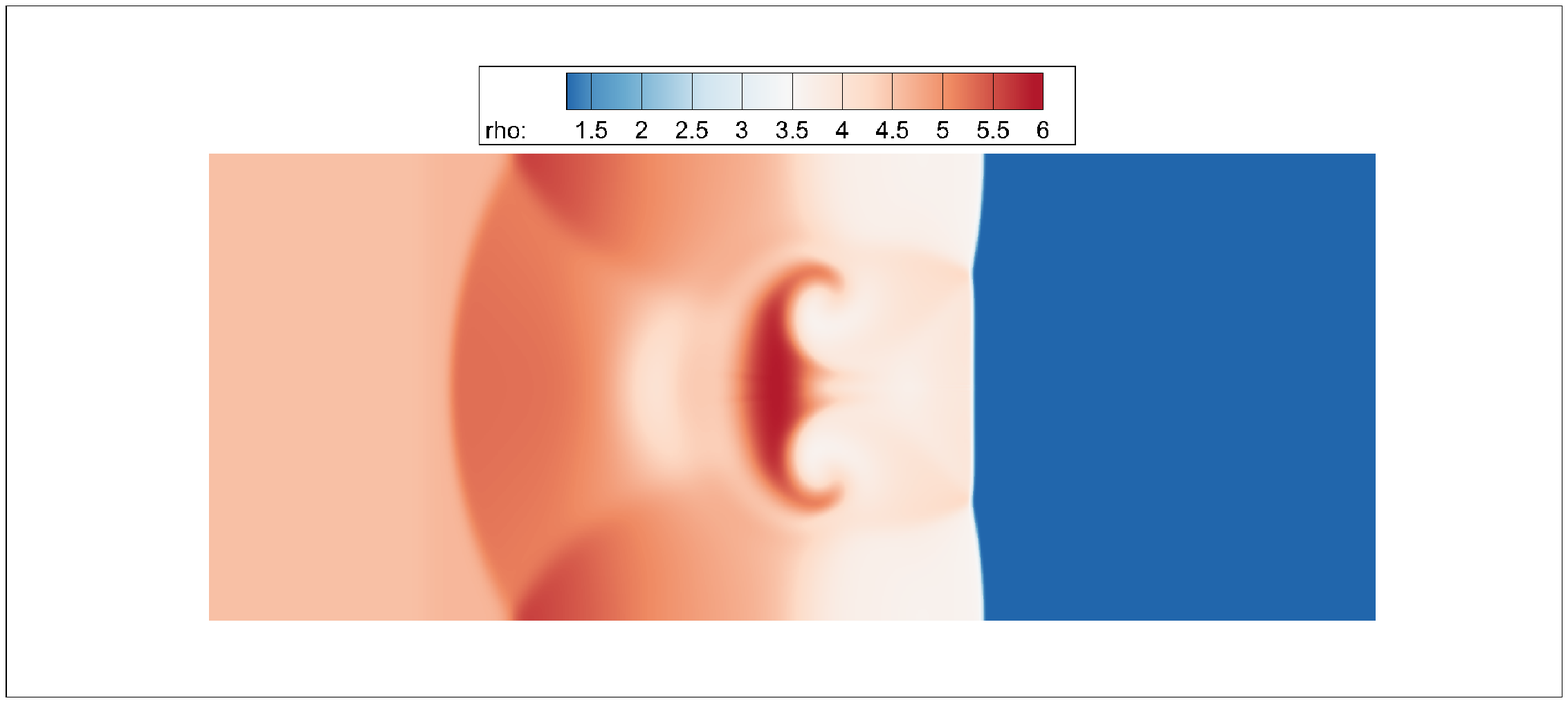}
  \includegraphics[width=7cm]{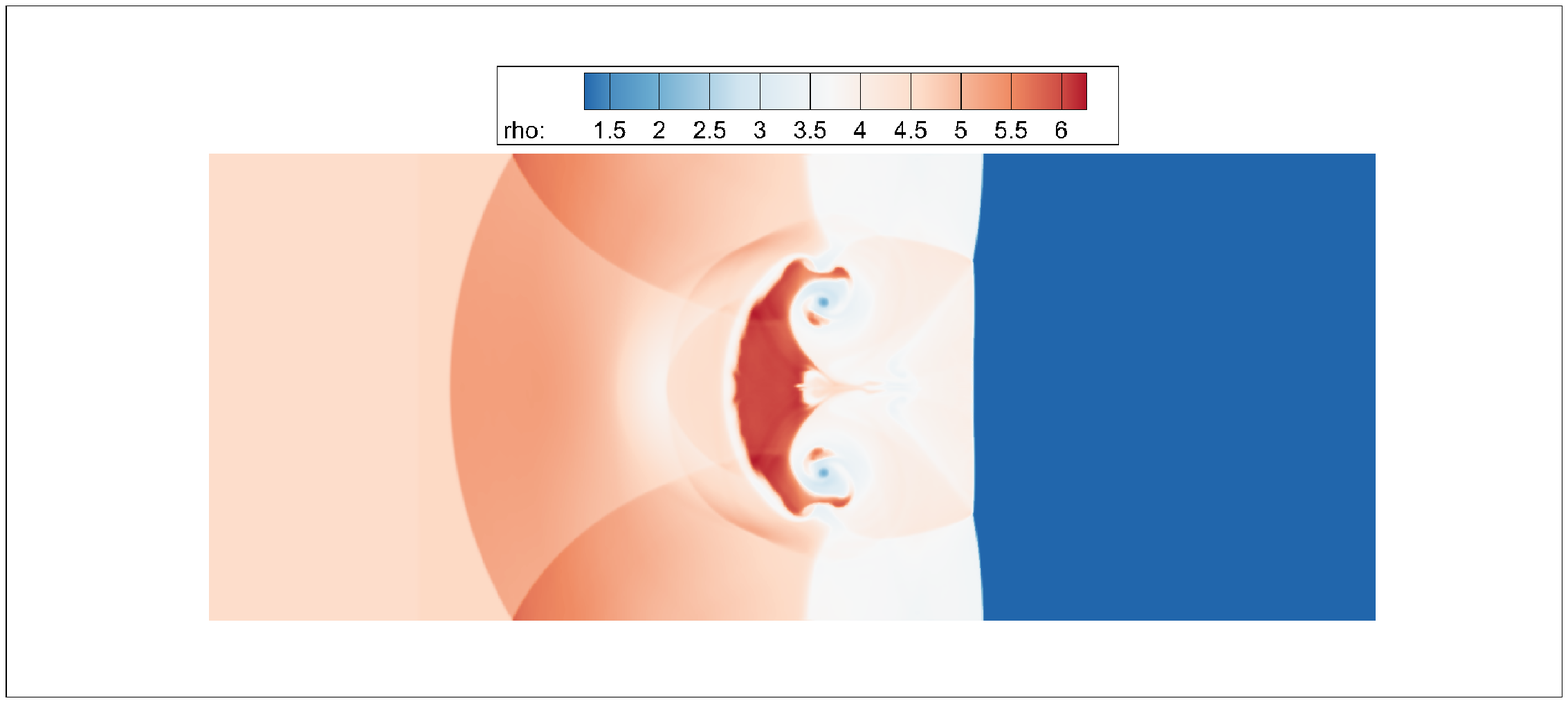}\\
  \includegraphics[width=7cm]{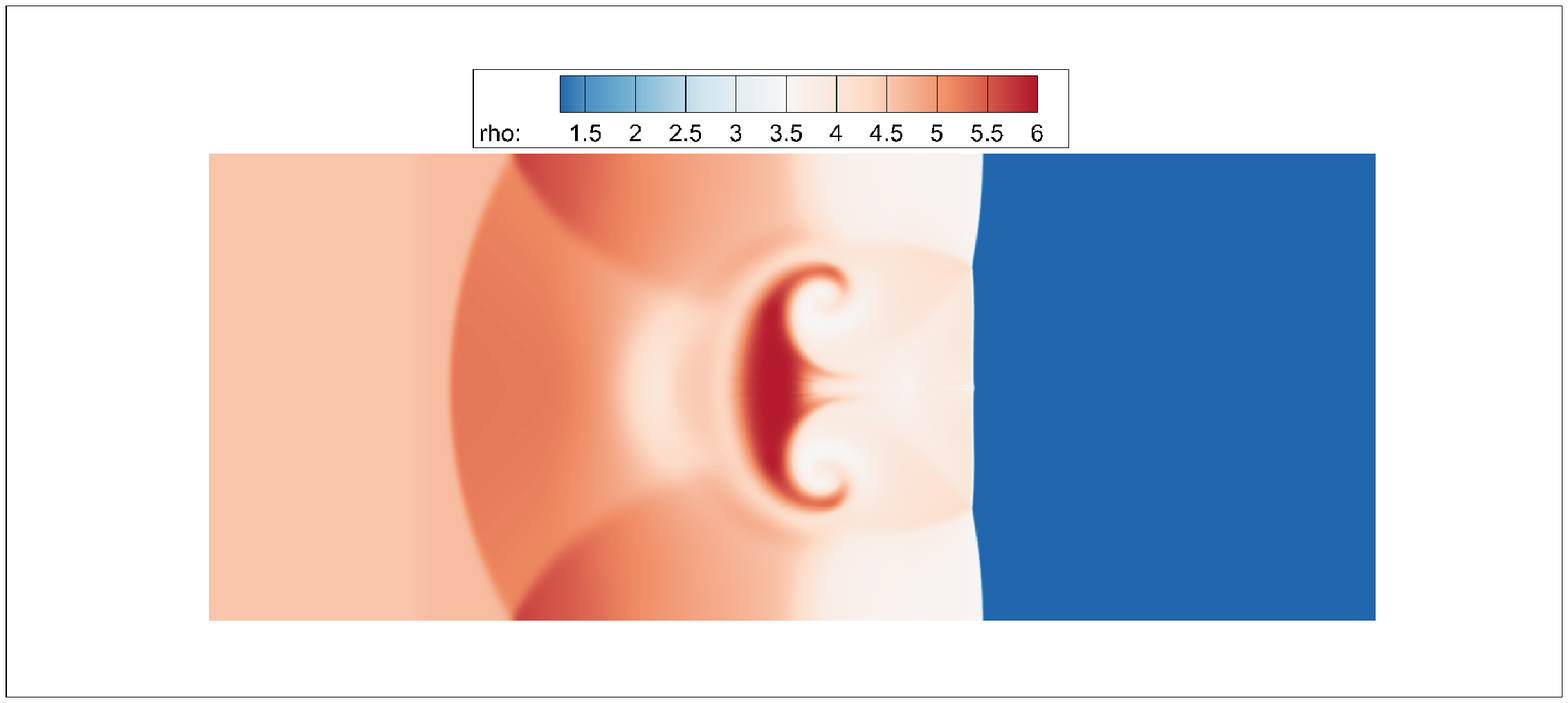}
    \includegraphics[width=7cm]{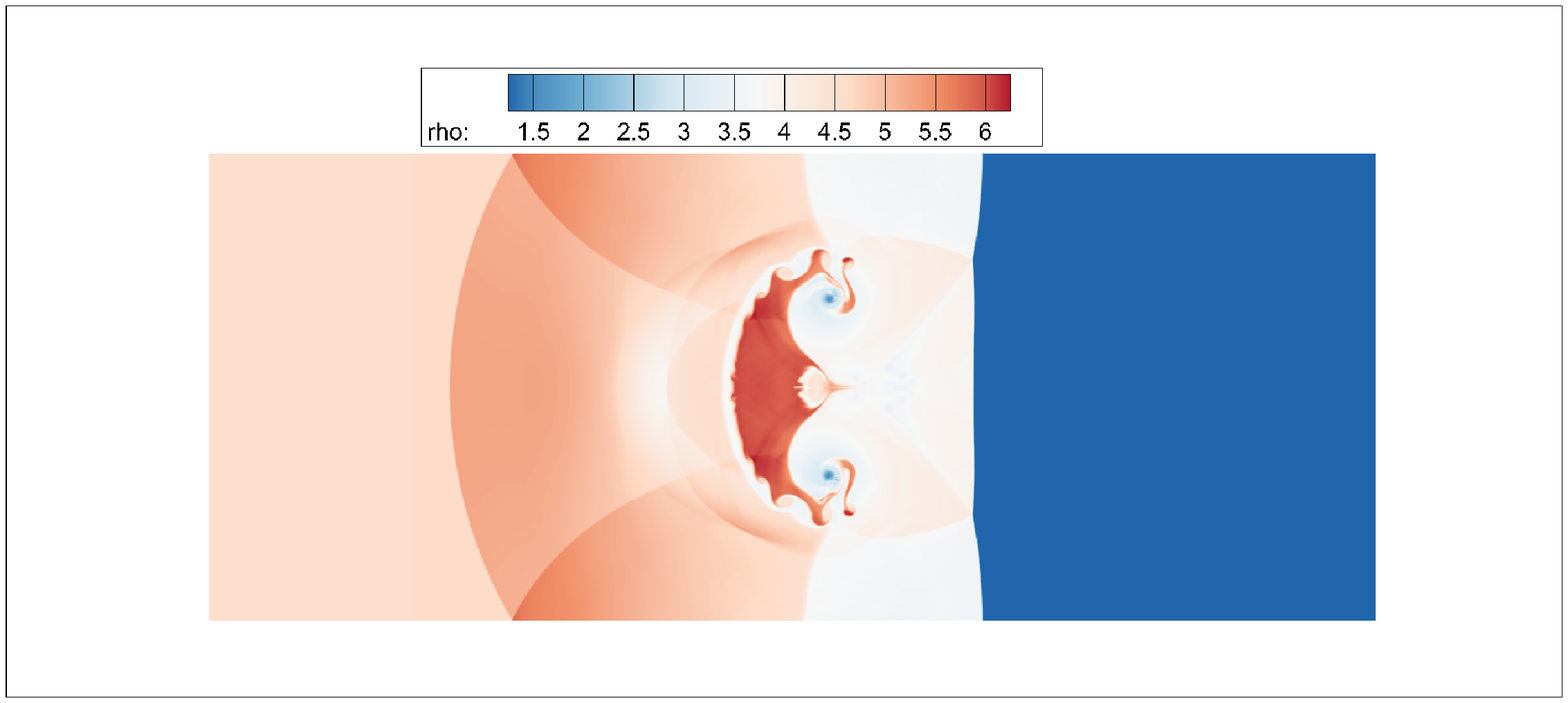}\\
  \caption{Shock-bubble interaction problem at $t = 70\mu s$: computed by the Godunov (left) and GRP (right) method with the stiffened gas approximation (Upper: $1200\times 400$ cells; Lower:$2400 \times 800$ cells). }\label{bubble-fig1}
\end{figure}
 \section{More remarks}

 In this paper we present Godunov type schemes for general EOS by using stiffened gas approximation and highlight the thermodynamic effects of the GRP scheme.
 Although this work is done just for a single material, the idea is of general significance.
 As for the use of the approximate GRP solver and its extension, some remarks are in order.
 \begin{enumerate}
 \item[(i)]  The original nonlinear GRP solver is useful for strong nonlinear waves, since thermodynamic effect becomes significant. The GRP solver with stiffened gas EOS approximation brings a new idea for engineering applications in order to make the scheme more robust and efficient while the thermodynamic effects are also embodied into the scheme.

 \item[(ii)] Following the same analogy, the GRP solver can be extend to multi-material fluid flows and the arbitrary Lagrangian-Euler (ALE) framework. For a specific EOS, the Gibbs relation has the corresponding implication which will be reflected to the corresponding numerical fluxes.

 \item[(iii)] For now, we only consider the EOS with the form \eqref{EOS}. It would be significant for engineering application to use the similar idea for solving fluid with more general EOSs, such as $p=p(\rho,T)$ for  the virial EOS and the van der Waals gas \cite{Plohr}  or the tabular forms.
  \end{enumerate}

\section*{Acknowledgment}  Jiequan Li is supported by NSFC (Nos. 11771054, 12072042,91852207), the Sino-German Research Group Project (No. GZ1465)  and Foundation of LCP.

\vspace{0.4cm}
\vspace{0.4cm}
\appendix

\section{The details of the stiffened gas approximate GRP solver}

In this appendix we assume a typical configuration (Rarefaction-Contact-Shock) and extract the detailed formulation of the two-material stiffened gas GRP solver.
As convention \cite{Ben-Artzi-1984,Ben-Artzi-2003,Ben2006}, the GRP solver  consists four parts, consistent with the associated Riemann solver: {\em (A) Resolution of rarefaction waves; (B) Tracking of shocks; (C) Bridge by contact discontinuity; and (D) Compute the time derivative of density$\rho$.}
More details can be found in   \cite{Li-Wang}.

\n {\bf (A)  Resolution of rarefaction waves. }

The key idea is to resolve the singularity at $(0,0)$. It consists of three steps:
\begin{enumerate}
\item[(I)] {\em Measurement of the expansion of the rarefaction wave.} Characteristic coordinates are applied to characterize the expansion of such a rarefaction wave. Let $\al(x,t)=C_1$ and $\beta(x,t)=C_2$ be the integral curves, respectively, of
  \begin{equation}
  \dfr{dx}{dt} =u+c,  \ \ \ \dfr{dx}{dt} =u-c. \notag
  \end{equation}
Denote $\Theta(\beta):=\frac{\pt t}{\pt \al}(0,\beta)$.  Then we have
\begin{equation}
\dfr{\Theta(\beta)}{\Theta(\beta_L)} = \left(\dfr{c(0,\beta)}{c_L}\right)^{\frac{1}{2\mu_L^2}}, \ \ \ \mu_L^2=\mu^2(\rho_L)=\frac{\gm(\rho_L)-1}{\gm(\rho_L)+1}.
\notag
\end{equation}
where $\beta_L=u_L-c_L$. $c(0,\beta_*)=c_{*L}$ when $\beta_*=u_*-c_{*L}$.

\item[(II)] {\em  The rate of entropy variation across curved rarefaction waves.}
 Across  the curved rarefaction wave associated with $u-c$, the entropy variation $T\pt S/\pt x(0,\beta)$ in the neighborhood of the singularity point has the change rate,
  \begin{equation}
  \dfr{T \pt S/\pt x(0,\beta)}{T_LS_{L}'}  =\left(\frac{c(0,\beta)}{c_L}\right)^{\frac{1}{\mu_L^2}+1}.
  \notag
  \end{equation}
Given the initial state in \eqref{data}, the initial variation of the entropy is  known thanks to the Gibbs relation \eqref{Gibbs},
  \begin{equation}
  T_L S_{L}' = e_{L}'-\dfr{p_L}{\rho_L^2} \rho_{L}'.
  \notag
  \end{equation}

 \vspace{0.2cm}

 \item[(III)] {\em The interaction  of kinematics and thermodynamics.} The entropy variation strongly affects the dynamics of kinematic variables.
  Apply \eqref{diagonal} and
return to the $(x,t)$-frame, then we have
\begin{equation}
a_L\dfr{D u}{D t}(0,\beta) +b_L \dfr{D p}{Dt}(0,\beta)  =d_L(0,\beta),
\notag
\end{equation}
where the total (material)  derivative $D/Dt =\pt/\pt t+u\pt/\pt x$, $\theta(\beta)=c(0,\beta)/c_L$ and
 \begin{align}
&a_L = 1, \ \ \ b_L = \dfr{1}{\rho(0,\beta) c(0,\beta)},\ \ \ \psi'_L=u'_L+\frac{1}{\rho_L c_L}p'_L+\frac{1}{c_L}T_LS'_L,\notag\\
&d_L=\left[\dfr{1+\mu_L^2}{1+2\mu_L^2}
\left(\theta(\beta)\right)^{1/(2\mu_L^2)}+\dfr{\mu_L^2}{1+2\mu_L^2}
\left(\theta(\beta)\right)^{(1+\mu_L^2)/\mu_L^2}\right]T_LS'_L-c_L\left(\theta(\beta)\right)^{1/(2\mu_L^2)}
\psi'_L.
\notag
\end{align}

\end{enumerate}

\n {\bf  (B) Tracking of the shock.}

The resolution of the right shock consists of the singularity tracking technique and the Lax-Wendroff strategy.
Let $x=x(t)$ be a shock with speed $\sg=x'(t)$ and separate two states $\bU_R(x,t)$ in the wave front and $\bU_*(x,t)$ in the wave back.
For the time being, denote $\bU_*(t) = \bU_*(x(t),t)$ and $\bU_R =\bU_R(x(t),t)$.
Using the Rankine-Hugoniot relations, the singularity is tracked by making differentiation along the shock trajectory $x=x(t)$. Denote
\begin{equation}
\dfr{D_\sg}{Dt }  = \dfr{\pt }{\pt t} +\sg\dfr{\pt }{\pt x}.\notag
\end{equation}
We specify to the shock associated with $u+c$.
The tracking of the right shock can be divided into two steps.
\begin{enumerate}
\item[(I)] {\em The interaction  of kinematics and thermodynamics.}
The Rankine-Hugoniot relations provide relations between all thermodynamic variables
\begin{align}
\rho_*=G(p_*;p_R,\rho_R)&=\rho_R\cdot\frac{p_*+\mu_R^2p_R+(1+\mu_R^2)p^R_{\infty}}{p_R+\mu_R^2p_*+(1+\mu_R^2)p^R_{\infty}},\label{RH2}
\end{align}
and interaction of kinematics and thermodynamics
\begin{align}
&u_*=u_R+\Phi(p_*;p_R,\rho_R)=u_R+(p^*-p_R)\Lambda^{\frac{1}{2}},\label{RH1}
\end{align}
where $\Lambda=(1-\mu_R^2)\tau_R/(p^*+\mu_R^2p_R+(1+\mu_R^2)p_{\infty}^R)$,
$p_{\infty}^R=p_{\infty}(\rho_R)$ and $\mu_R^2=\mu^2(\rho_R)=\frac{\gm(\rho_R)-1}{\gm(\rho_R)+1}$.

\item[(II)] {\em  Tracking    the shock trajectory $x=x(t)$.}
Take the differentiation along the shock trajectory $x=x(t)$,
\begin{equation}
\dfr{D_\sg u_*}{Dt} = \dfr{D_\sg  u_R}{Dt}  + \dfr{\pt \Phi}{\pt p_*} \dfr{D_\sg p_*}{Dt} +  \dfr{\pt \Phi}{\pt  \rho_R} \dfr{D_\sg  \rho_R}{Dt} + \dfr{\pt \Phi}{\pt p_R} \dfr{D_\sg  p_R}{Dt},\notag
\end{equation}
 and the limit along the shock trajectory $x=x(t)$, we have
\begin{equation}
a_R\left(\dfr{D u}{Dt}\right)_* +b_R\left(\dfr{D p}{Dt}\right)_* =d_R,\notag
\end{equation}
where the coefficients $a_R$, $b_R$ and $d_R$ are given in terms of the intermediate state $\bU_*$ and the initial data from the right,
\begin{equation}
\begin{array}{l}
a_R=1+\rho_{*R}(\sg-u_*)\dfr{\pt \Phi}{\pt p_*}(p_*;p_R,\rho_R),\\[3mm]
b_R=-\left[\dfr{\sg-u_*}{\rho_{*R}c_{*R}^2}+\dfr{\pt \Phi}{\pt p_*}(p_*;p_R,\rho_R)\right],\\[3mm]
d_R=L_p^R \cdot  p_R'+L_u^R \cdot u_R'+L_{\rho}^R \cdot \rho_R', \\\\
\end{array}
\notag
\end{equation}
where
\begin{equation}
\begin{array}{l}
 L_p^R= -\dfr{1}{\rho_R}+(\sg-u_R)\dfr{\pt \Phi}{\pt
 p_R}(p_*;p_R,\rho_R),\\[3mm]
 L_u^R=
\sg-u_R -\rho_R c_R^2\dfr{\pt \Phi}{\pt p_R}(p_*;p_R,\rho_R)-\rho_R \dfr{\pt
\Phi}{\pt \rho_R}(p_*;p_R,\rho_R),\\[3mm]
L_{\rho}^R= (\sg-u_R)\dfr{\pt
\Phi}{\pt \rho_R}(p_*;p_R,\rho_R).
\end{array}
\notag
\end{equation}
and
\begin{equation}
\begin{array}{l}
 \dfr{\pt \Phi}{\pt p_*}=\frac{\sqrt{\Lambda}}{2}\frac{p_*+(1+2\mu_{R}^2)p_R+2(1+\mu_R^2)p^R_{\infty}}{(p^*+\mu_R^2p_R+(1+\mu_R^2)p^R_{\infty})^2},\\[3mm]
 \dfr{\pt \Phi}{\pt p_R}=-\frac{\sqrt{\Lambda}}{2}[\frac{(\mu_R^2+2)p_*+\mu_R^2p_R+2(1+\mu_R^2)p^R_{\infty}}{p^*+\mu_R^2p_R+(1+\mu_R^2)p^R_{\infty}}],\\[3mm]
 \dfr{\pt \Phi}{\pt \rho_R}=-\frac{p_*-p_R}{2\rho_R}\sqrt{\Lambda}.
\end{array}
\notag
\end{equation}
\end{enumerate}

\vspace{0.2cm}
\n{\bf (C) Bridge by contact discontinuity. }

Thanks to the continuity property of $u$ and $p$ across the contact discontinuity, their material derivatives  keep continuous, respectively. It turns out that the contact discontinuity bridges two sub-regions between the shock and the rarefaction waves,
\beqn
\left\{
\begin{array}{l}
a_L\left(\dfr{D u}{Dt}\right)_* +b_L\left(\dfr{D p}{Dt}\right)_* =d_L,\\[3mm]
a_R\left(\dfr{D u}{Dt}\right)_* +b_R\left(\dfr{D p}{Dt}\right)_* =d_R.
\end{array}
\right.
\notag
\eeqn
Solving this algebraic system provides $\left(\dfr{D u}{Dt}\right)_*$ and $\left(\dfr{D p}{Dt}\right)_*$. Then we obtain
\begin{align}
&\left(\dfr{\pt u}{\pt t}\right)_*=\left(\dfr{D u}{Dt}\right)_*+\frac{u_*}{\rho_*c_*^2}\left(\dfr{D p}{Dt}\right)_*,\notag\\
&\left(\dfr{\pt p}{\pt t}\right)_*=\left(\dfr{D p}{Dt}\right)_*+\rho_*u_*\left(\dfr{D u}{Dt}\right)_*.\notag
\end{align}

\n {\bf  (D) Derive the time derivative for the density $\rho$}

If $u_*>0$, we return to the rarefaction wave by using $S_t=-uS_x$ to derive
\begin{align}
(\frac{\pt \rho}{\pt t})_*&=c_{*L}^{-2}((\pt p/\pt  t)_*+(\gm_L-1)\rho_{*L}u_*T_LS_L^{\prime}(\frac{c_{*L}}{c_L})^{1+\frac{1}{\mu_L^2}}).\notag
\end{align}
where $\gm_L=\gm(\rho_L)$.

If $u_*<0$, the limiting value $(\frac{\pt \rho}{\pt t})_*$ is derived by the Rankine-Hugoniot relation
\begin{align}
&\tilde{G}(\rho_{*R},p_*;p_R,\rho_R)=\frac{p_*+\gm_R p_{\infty}^R}{(\gm_R-1)\rho_{*R}}-\frac{p_R+\gm_R p_{\infty}^R}{(\gm_R-1)\rho_R}-\frac 12(p_R+p_*)(\frac{1}{\rho_R}-\frac{1}{\rho_{*R}})=0,\notag\\
&(u_*-\sigma)H_1(\frac{\pt \rho}{\pt t})_*+(\frac{\sigma}{c_{*R}^2}H_1+u_*H_2)(\frac{Dp}{Dt})_*+H_2(u_*-\sigma)\rho_{*R}u_*(\frac{Du}{Dt})_*=\notag\\
&u_*((u_R-\sigma)H_3\rho'_R+(\rho_RH_3+\rho_R c_R^2H_4)u'_R+(u_R-\sigma)H_4p'_R).\notag
\end{align}
where $\gm_R=\gm(\rho_R)$, $p_{\infty}^R=p_{\infty}(\rho_R)$ and
\begin{align}
&H_1=\frac{\pt\tilde{G}}{\pt\rho_*}=-\frac{1}{2}(\frac{p_*}{\rho_{*R}^2\mu_R^2}+\frac{p_R}{\rho_{*R}^2}+(1+\mu_R^2)\frac{p_{\infty}^R}{\rho_{*R}^2}),\notag\\
&H_2=\frac{\pt\tilde{G}}{\pt p_*} = \frac 12(\frac{1}{\rho_{*R}\mu_R^2}-\frac{1}{\rho_R}),\notag\\
&H_3=\frac{\pt\tilde{G}}{\pt\rho}=\frac{1}{2}(\frac{p_R}{\rho_R^2\mu_R^2}+\frac{p_*}{\rho_R^2}+(1+\mu_R^2)\frac{p_{\infty}^R}{\rho_R^2}),\notag\\
&H_4=\frac{\pt\tilde{G}}{\pt p} =-\frac{1}{2}(\frac{1}{\rho_R\mu_R^2}-\frac{1}{\rho_{*R}}).\notag
\end{align}


\begin{thebibliography}{10}

\bibliographystyle{plain}

\bibitem{Bai}
J. S. Bai, Z. J. Zhang, P. Li, M. Zhong, Numerical method of the Riemann problem for two-dimensional multi-fluid flows with general equation of state,
{\em Chin. Phys.}, 15 (2006), 22--34.

\bibitem{banks}
J. W. Banks, On exact conservation for the Euler equations with complex equations of state, {\em Commun. Comput. Phys.}, 8 (2010), 995--1015.

\bibitem{Ben-Artzi-1984}
M. Ben-Artzi and J. Falcovitz, A second-order Godunov-type scheme for compressible fluid dynamics,  {\em J. Comput. Phys.}, 55 (1984),  1--32.

\bibitem{Ben-Artzi-2003}
M. Ben-Artzi and J. Falcovitz, Generalized Riemann problems in computational gas dynamics, {\em Cambridge University Press}, 2003.

\bibitem{Ben2006}
M. Ben-Artzi, J. Li and G. Warnecke, A direct Eulerian GRP scheme for compressible fluid flows, { \em J. Comput. Phys.}, 218 (2006), 19--34.

\bibitem{Ben2007}
M. Ben-Artzi and J. Li, Hyperbolic balance laws: Riemann invariants and the generalized Riemann problem, {\em Numer. Math.}, 106 (3) (2007), 369--425.





\bibitem{Godunov}
S. K. Godunov, A finite difference method for the numerical  computation and disontinuous solutions of the equations of fluid dynamics, {\em Mat. Sb.}, 47 (1959), 271--295.

\bibitem{Godunov1979}
S. K. Godunov, A. Zabrodine, M. Ivanov, A. Kraiko, and G. Prokopov, R\'esolution num\'erique des probl\`emes
multidimensionnels de la dynamique des gaz (Editions Mir, Moscow, 1979).


\bibitem{Gottlieb}
J. Gottlieb, C. Groth, Assessment of Riemann solvers for unsteady one-dimensional inviscid flows of perfect gas, {\em J. Comput.
Phys.}, 78 (1988) 437--458.


\bibitem{Kamm}
J. R. Kamm, An exact, compressible one-dimensional Riemann solver for general, convex equation of states, {\em LA-UR-15-21616}, 2015.

\bibitem{Lee}
B. J. Lee, E. F. Toro, C. E. Castro, and N. Nikiforakis, Adaptive Osher-type scheme
for the Euler equations with highly nonlinear equation of states, {\em J. Comput. Phys.}, 246 (2013), 165--183.



\bibitem{Li-Wang}
J. Li and Y. Wang, Thermodynamical Effects and High Resolution Methods for Compressible Fluid Flows, {\em J. Comput. Phys.}, 343 (2017), 340--354.

\bibitem{Plohr}
R. Menikoff ad B. J.  Plohr, The Riemann problem for fluid flow of real materials, {\em Rev. Modern}, 61 (1) (1989), 75--130.




\bibitem{Quar-2003}
L. Quartapelle, L.  Castelletti, A. Guardone and G. Quaranta, Solution of the Riemann problem of classical gasdynamics,
{\em J. Comput. Phys.},  190 (1) (2003), 118--140.


\bibitem{Saurel}
R. Saurel, R. Abgrall, A Multiphase Godunov Method for Compressible Multifluid and Multiphase Flows, {\em J. Comput. Phys.}, 150 (1999), 425--467.

\bibitem{Saurel2007}
R. Saurel, E. Franquet, E. Daniel, O. Le Metayer, A relaxation-projection method for compressible flows. I. The numerical equation of state for the Euler equations, {\em J. Comput. Phys.}, 223 (2007), 822--845.

\bibitem{Shyue}
K.-M. Shyue, A fluid-mixture type algorithm for compressible multicomponent flow with Mie-Gr{\"{u}}neisen equation of state, {\em J. Comput. Phys.}, 171 (2001), 678--707.

%

\bibitem{Tang-Liu}
H. Z. Tang and T. G. Liu, A note on the conservative schemes for the Euler equations, {\em J. Comp. Phys.}, 218 (1) (2006), 451--459.

\bibitem{Thompson}
P. A. Thompson, Compressible-Fluid Dynamics, {\em Hamilton Press}, 1988.

\bibitem{Toro}
E. F. Toro, Riemann solvers and numerical methods for fluid dynamics. A practical introduction, {\em Springer--Verlag, Berlin}, 1997.



\end{thebibliography}
\end{document}